\documentclass[]{interact}

\usepackage{amsmath}
\usepackage{graphicx}     
\usepackage{hyperref} 
\usepackage{url}
\usepackage{amsfonts} 
\usepackage{multicol}
\usepackage{float}
\usepackage{subfig}
\usepackage{array}
\usepackage{booktabs}

\theoremstyle{plain}
\newtheorem{theorem}{Theorem}[section]
\newtheorem{lemma}[theorem]{Lemma}

\newtheorem{proposition}[theorem]{Proposition}

\theoremstyle{definition}

\theoremstyle{remark}
\newtheorem{remark}{Remark}

\allowdisplaybreaks

\makeatletter
\@addtoreset{footnote}{page}
\makeatother

\begin{document}
\title{\LARGE A random journey through the math of gambling.}
\markright{An unexpected random journey.}

\author{\Large Paolo Bartesaghi\\               
\scriptsize Department of Statistics and Quantitative Methods\\
University of Milano - Bicocca\\
via Bicocca degli Arcimboldi 1, 20100 Milano (Italy)\\    
paolo.bartesaghi@unimib.it\\
paolo.bartesaghi73@gmail.com\\
ORCID ID 0000-0001-6539-2330}                      

\maketitle

\begin{abstract}
The laws of chance are often subtle and deceptive. This is why games of chance work. People are convinced that they obey seemingly intuitive laws, while the underlying mathematical structure reveals a different and more complex reality. This article is a brief and rigorous journey through the implications that the mathematical laws governing stochastic processes have on gambling. It addresses a specific process, the random walk, and analyze some instances of fair and unfair games by highlighting the fallacy of many of our intuitions and beliefs. The paper gradually moves from the analysis of the random walk properties to a comprehensive description of the ruin problem. The introduction of the idea of transient and persistent states concludes the discussion.  Much emphasis is placed on concrete examples and on the numerical values, in particular of the involved probabilities, and the interpretation of the results is always more central than the demonstrative technical details, which are nevertheless available to the reader.
\end{abstract}



\section{Introduction}
\label{section1}
The simplest game of chance between two players is to flip a coin so that one wins on heads and the other wins on tails. Suppose the two players repeat the toss many times, and each time one loses, he or she gives the other a dollar. Of course, a player is considered a winner at some stage of the game if he has won more than he has lost. Now consider the question: What fraction of the total game time do we imagine one player will spend in the role of the winner, and what fraction in the role of the loser?

Most people will probably respond that, being tossing a coin a '50\%-50\%' game, each player will spend about half the time leading and half the time losing. Indeed we are more or less convinced that we will spend part of the time as winners and part of the time as losers, and that this distribution will be somehow balanced. We also know that we cannot expect a $50\%$ and $50\%$ net split and that there will be some kind of fluctuation. Maybe we can expect to be in the lead $55\%$ of the time and our opponent $45\%$ of the time, or $48\%$ against a $52\%$. In short, although we do not know what the exact split will be, we think it is most likely to be around half the playing time.

In the same way we are led to think that, if we start to loose, sooner or later we will break even and regain what we lost. But how many times can we afford to lose and come back to zero? Are we sure that this can happen as often as we would like?

Furthermore, what if the coin is rigged and it is no longer a fair game? how should we adjust our expectations? in a sense, in this case, it is easier for us to imagine some kind of asymmetric behavior that makes us in the game more favorable or unfavorable depending on the point of view.

Let us mention here only one example that will be discussed in the text. Suppose our strategy is to quit the game when our total winnings reach a certain value $A$, or, alternatively, when we completely deplete the budget $B$ that we decide to allocate to this game. We will see that if even the probability of winning at each step is only 5\% below the 50\% which is typical of a fair game, and the probability of losing is 5\% higher, then the probability of winning $A$ before losing all of $B$ is substantially less than the probability of losing all of $B$ before winning $A$. Suppose $A$ and $B$ are 3 dollars. Then the probability of losing all our small 3 dollars budget before winning other 3 dollars would be $65\%$, as opposed to $35\%$ of winning all 3 dollars before losing them. If it were 10 dollars, this probability would rise to $89\%$, and the probability of losing would go up even more for higher bets.

This article will attempt to shed light on these and some other similar questions and to show that issues are more subtle than they appear on the surface. It is designed for anyone who wants to get a clearer idea of the mechanisms that govern gambling and games of chance in their mathematical implications, which, as paradoxical, also become psychological. We begin with a seemingly unrelated problem, counting votes in an election between two candidates, and then describe the properties of the random walk in the case of a fair game. We then introduce the ruin problem and compute the probabilities mentioned above in both the fair and unfair cases. Finally, we conclude with some remarks on the persistence of the random walker's return to the origin.

Reading requires a certain amount of patience. The text deliberately aims for maximum accessibility combined with minimum mathematical rigor, so that there are no logical or narrative leaps. The mathematics used is essentially combinatorics and elementary probability. The reader will not find any 'it is clearly seen that' or 'it is left to the reader to demonstrate immediately' except in one or two places at most where things are really obvious.
However, the reader will not find any proofs in the text. It was simply felt that all proofs should be placed in the appendix, so that the text is direct and aimed at understanding the ideas rather than analyzing the demonstrative details. The reader who wants to go deeper is thus directed to the appendix. Only one proof is kept in the text since it is very instructive.

\section{The Ballot Problem}
\label{section2}

Let us imagine that, during the ballot for the election to an institutional office, there are two candidates $P$ and $Q$ and that the ballot proceeds by examining the voting records one by one. Of course, during the counting process, one candidate may be first in the lead and then the other may move ahead, and this switch may occur several times. The question we want to focus on is: What is the probability that, throughout the counting, there have always been more votes for $P$ than for $Q$?

A good way to list the votes for each candidate, at each step $n$ of the counting process, is as follows. We can mark on a blackboard a $+1$ when the vote is assigned to $P$ and a $-1$ when the vote is assigned to $Q$. After $n$ steps, we have a list where we filled in imaginary cells with a sequence of $+1$ or $-1$. The value assigned to each cell is completely random. We call the empty cells \textit{random variables} and denote them by $X_i$. We have $n$ random variables that get the value $+1$ if the vote is assigned to $P$ and the value $-1$ if the vote is assigned to $Q$:
\begin{equation}
	X_i =\left\{ 
	\begin{array}{cc}
		+1 & {\rm if}\ \ {\rm vote\ assigned\ to}\ P \\ 
		-1 & {\rm if}\ \ {\rm vote\ assigned\ to}\ Q
	\end{array}
	\right.
\end{equation}
We are left with a sequence $\{ X_1, \dots, X_n \}$ given by the votes we get from the ballot at each step. Suppose that in the sequence we have $p$ plus one ($+1$) and $q$ minus one ($-1$), and that $n$ is the last step, or, equivalently, that we are interested only in what happened before step $n$. Of course, the partial sum
\begin{equation}
	S_k=\sum_{i=1}^{k}X_i
\end{equation}
is the number of votes by which $P$ leads or trails just after the $k-$vote is cast. Clearly we have 
$S_n=p\cdot(+1)+q\cdot (-1)=p-q$ and $S_k-S_{k-1}=X_k=\pm 1$ with $S_0=0$, for any $k=1, \dots, n$. Conversely, an arrangement $\{S_1, \dots, S_n\}$ of integers satisfying $S_k-S_{k-1}=\pm 1$ can represent a potential voting record.

\begin{figure}[H]
	\includegraphics[width=0.75\linewidth]{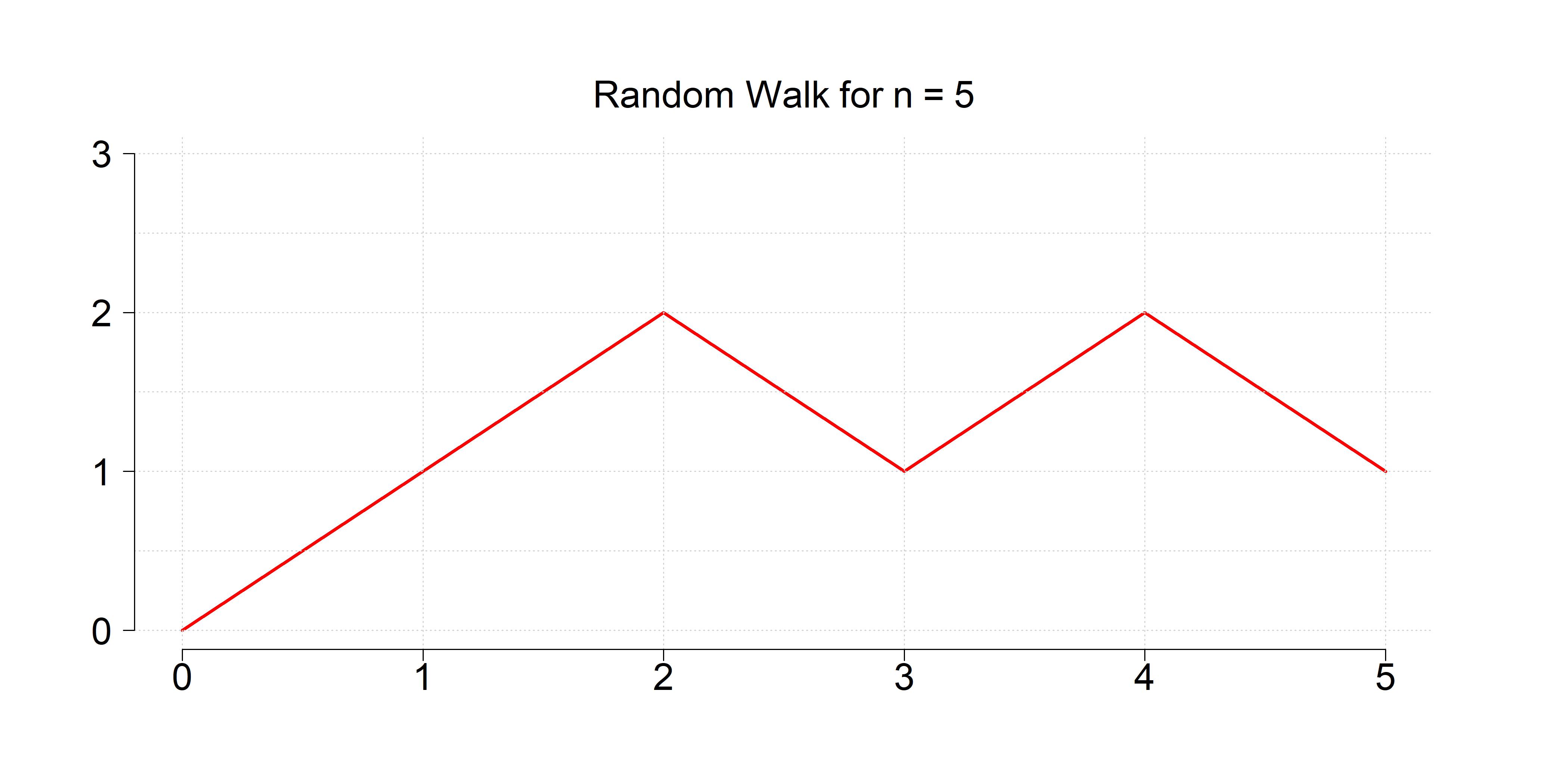}
	\centering
	\caption{A simple example of possible outcome of the ballot counting.}
	\label{fig1}      
\end{figure}

Now the idea is to represent such an arrangement by a polygonal line whose $k-$th vertex has ordinate $S_k$ and whose $k-$th side has slope $X_k$. This line is called {\it path}.
Consider the simple example illustrated in Fig. \ref{fig1}: after the counting of five votes, that is $n=5$, we have $S_n=1$. Therefore $p+q=n=5$ and $p-q=S_n=1$, which means $p=3$ and $q=2$. The figure represents the path of the counting process, or, in other terms, of the random process.
\vskip -0.5cm

More in general, a path $\{S_1, \dots, S_n\}$ from the origin to the point $(x,y)=(n,S_n)$ is a polygonal line whose vertices are $(k,S_k)$ with $S_k-S_{k-1}=X_k$ and $y=S_n$.

Since at each of the $n$ steps we have two possible directions, \textit{up} or \textit{down}, there exist $2^n$ different paths coming out from the origin in $n$ steps to the end point $(n,y)$, for an arbitrary, that is not fixed, integer $y$.
If, on the other hand, the value of $y$ is constrained to be of the type $p-q$ when $x$ is equal to $p+q$, the number of possible paths is significantly reduced. In fact, if exactly $p$ among the $X_k$ are $+1$ and $q$ are $-1$, so that $x=n=p+q$ and $y=S_n=p-q$, then the $p$ places in which we can put the $+1$, can be chosen from the total $n=p+q$ available places in $\binom{n}{p}=\binom{n}{q}$ different ways. Therefore, this number is equal to
\begin{equation}
	N_{x,y}=\binom{p+q}{p}=\binom{p+q}{q}=\binom{x}{\frac{x+y}{2}}=\binom{x}{\frac{x-y}{2}},
\end{equation}
and $N_{x,y}=0$ if $x$ and $y$ are not of the form $x=p+q$ and $y=p-q$. This means that there are exactly $N_{x,y}$ different paths from the origin to the point $(x,y)=(p+q,p-q)$. Now, the condition that there have always been more votes for $P$ than for $Q$ is equivalent to the fact that $S_1>0$, $S_2>0$, $\dots$, $S_{n-1}>0$ and $S_{n}=y>0$, i.e. that the path is always positive. 
The fraction of all the $N_{x,y}$ paths, from the origin to point $(x,y)=(p+q,p-q)$, that are always positive is quantified in a well-known result, named after J. L. F. Bertrand (1887, see \cite{Feller}), also called the Ballot Theorem:
\begin{proposition}
	If $x>0$ and $y>0$, i.e. $n>0$ and $S_n>0$, the expected number of paths from the origin to the point $(x,y)$ that are always positive, i.e. of the kind $\{S_1,S_2, \dots, S_n=y\}$ with $S_i>0, \forall i=1, \dots n$ equals
	\begin{equation}
		\frac{p-q}{p+q}N_{x,y}=\frac{y}{x}N_{x,y}
	\end{equation}
	\label{theorem1}
\end{proposition}
\vskip -0.5cm
This lemma provides the probability that, if $P$ is leading at step $n$, then \textit{it has always led at all previous steps}. The complete proof of the Proposition \ref{theorem1} can be found in the appendix. Here we make a few remarks. First, let us consider again the example in Fig. \ref{fig1}. We said that $x=n=5$, and $y=S_n=1$, so that $p=3$ and $q=2$.
Then there are $N_{x,y}=\binom{5}{3}=10$ possible paths from $(0,0)$ to $(5,1)$, but only $\frac{y}{x}N_{x,y}=\frac{1}{5}\cdot 10=2$ paths are always positive and they are the two sequences $(0,1,2,1,2,1)$ and $(0,1,2,3,2,1)$. Note that $\frac{1}{5}$ is the slope of the straight line from $(0,0)$ to $(5,1)$ and it is the probability that $P$ \textit{it has always led at the previous steps} if he or she is leading at step $n=5$. In general, the total number of possible paths from the origin to the point $(x,y)$ is multiplied by the mean slope of the polygonal line, that is the slope of the straight line from the origin to that point.

\hfill

\begin{figure}[H]
	\includegraphics[width=0.80\linewidth]{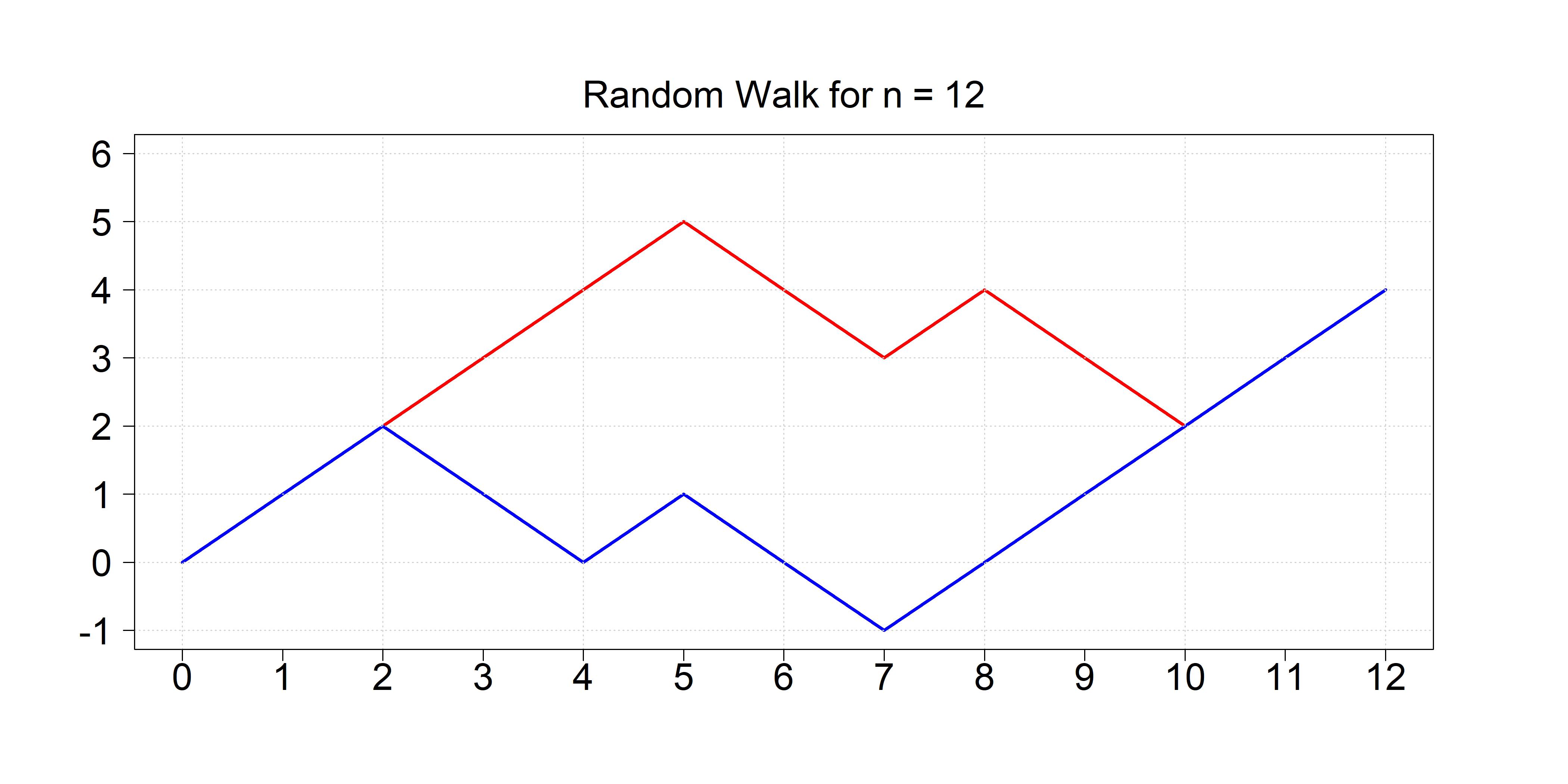}
	\centering
	\caption{A simple example of reversed path. Blue line is the path $S$ and red line (when non overlapping with $S$) is the reversed path $S^{\star}$.}
	\label{fig2}      
\end{figure}

{\bf A reverse journey.} Proposition \ref{theorem1} can be re-formulated in an interesting way in terms of {\it reversed path}. It is in a sense a dual formulation of that statement. First of all, we define {\it reversed path} the path obtained by reversing the order of the $X_i$'s. In the example in Fig. \ref{fig2}, the actual sequence of random variables in the blue path $S$ is $(+1,+1,-1,-1,+1,-1,-1,+1,+1,+1,+1,+1)$. The reversed path $S^{\ast}$ is given by the sequence $(+1,+1,+1,+1,+1,-1,-1,+1,-1,-1,+1,+1)$.

The reversed path is then described by the partial sums of the $X_i$'s in the reversed order, that is
\begin{equation*}
	\begin{split}
		S_{1}^{\ast} & = X_n = S_n-S_{n-1}\\
		S_{2}^{\ast} & = X_n+X_{n-1} = S_n-S_{n-2} \\
		S_{3}^{\ast} & = X_n+X_{n-1}+X_{n-2} = S_n-S_{n-3} \\
		& \ \ \vdots \\
		S_{n}^{\ast} & = X_n+X_{n-1}+X_{n-2}+ \dots +X_1 = S_n \\
	\end{split}
\end{equation*}

As can be seen in Fig. \ref{fig2}, the two paths are congruent, join the same endpoints and are obtained from each other by a rotation of $180^{\circ}$. Now proposition \ref{theorem1} can be re-formulated in this way. If $x,y>0$, the number of reversed paths $\{ S_{1}^{\ast}, \dots , S_{n}^{\ast} \}$ joining the point $(0,0)$ to $(x,y)$ such that $S_{i}^{\ast}>0$ for $i=1,2,\dots ,n$ is equal to $\frac{y}{x}N_{x,y}$. Note that asking for $S_{i}^{\ast}>0$ is equivalent to asking for $S_{n}-S_{n-i}>0$, i.e. $S_{n}>S_{n-i}$, for $i=1,2,\dots ,n$. Explicitly, this means $S_{1}<S_{n}$, $S_{2}<S_{n}$, $\dots$, $S_{n-1}<S_{n}$. Geometrically speaking, we can say that the original version of the theorem is concerned with paths whose left endpoint is the lowest vertex, whereas the dual version of the theorem is concerned with paths whose right endpoint is the highest.

{\bf Back to the origin.} Before concluding this section, we want to make some remarks about the special case where the endpoint is $y=0$, i.e. where the two candidates receive the same number of votes $p=q$ after $n$ rounds of voting. This is equivalent to a path that returns to the origin after its random walk above or below the axis. Since it is impossible to return to the $x$-axis in an odd number of steps, or equivalently, we cannot have an odd vertex on the $x$-axis, from now on we will denote this point as $(2n,0)$ and we will be interested in paths that connect $(0,0)$ to a point $(2n,0)$ on the $x$-axis. With this choice $p+q=2n$ and $p=q=n$. Of course, the number of all these paths is $N_{2n,0}=\binom{2n}{n}$. Let us first define $L_{2n}$ the following number
\begin{equation}
L_{2n}=\frac{1}{n+1}N_{2n,0}=\frac{1}{n+1}\binom{2n}{n}
\end{equation}

\begin{proposition}
The following statements hold
\begin{itemize}
\item[a)] The number of paths such that $S_{1}\geq 0$, $S_{2}\geq 0$, $\dots$, $S_{2n-1}\geq 0$ and $S_{2n}= 0$ is equal to $L_{2n}$; 
\item[b)] The number of paths such that $S_{1}> 0$, $S_{2}> 0$, $\dots$, $S_{2n-1}> 0$ and $S_{2n}= 0$ is equal to $L_{2n-2}$.
\end{itemize}
\label{theorem2}
\end{proposition}

Now the question is: Is it truly as easy as it might seem to return to the origin? What is the likelihood of this happening after a candidate has been in the lead for the entire ballot count? We will see that such an event is indeed extremely unlikely. To show this, we move from the problem of a runoff between two candidates to a more abstract game that reproduces the same conditions.

\section{Tossing coins}
\label{section3}
Let us now focus on a fair game, like flipping a coin to get heads or tails. Imagine that two friends $P$ and $Q$ decide to play the following game for a given period, let's say $N=20$ days. Every day they toss a coin. If it comes up heads, $P$ wins and $Q$ loses; if it comes up tails, $P$ loses and $Q$ wins. The value of the random variable is then assigned day by day, and each of the two friends constructs his own random walk.

Indeed, compared to the votes cast for candidates $P$ and $Q$, in this case we can make a \textit{a priori} assumption about the probability of getting one or the other result. If we call $p$ and $q$ the two probabilities of getting head or tail, they will be the same: $p=q=\frac{1}{2}$. Suppose we decide that a particular random variable $X_{i}$ takes on the value $+1$ if heads comes out and the value $-1$ if tails comes out, then we can say that the probabilities of the random variable taking on either value are given by
\begin{equation}
{\cal P}(X_i=1)=\frac{1}{2}\qquad {\cal P}(X_i=-1)=\frac{1}{2}
\end{equation}

Therefore, we can represent any possible outcome of the $N$ successive tosses of a coin-and so any possible path followed by the two players in the N days of their game-by a path of $N$ sides starting at the origin. If we get heads, $P$ goes up and $Q$ down; if we get tails, $P$ goes down and $Q$ up. Conversely, each such path can be seen as representing the outcome of $N$ tosses of a coin. We have already said that the total number of possible paths of this type is $2^N$. The set of all these paths is called \textit{sample space}. The sample space is the collection $\Omega$ of the $2^N$ paths $\{ S_{1}, \dots , S_{N} \}$ starting at the origin. Since we have no reason to think that one of the paths is preferred over the others, we can attribute probability $2^{-N}$ to each one. 

An event such as $E=\{ {\rm two\ heads\ at\ the\ first\ two\ trials} \}$ must be interpreted as the aggregate of all sequences starting with $S_{1}=1$ and $S_{2}=1$. There are $2^{N-2}$ such sequences, since the first two steps are predetermined. The probability of the event $E$ is then ${\cal P}(E)=2^{N-2}/2^{N}=2^{-2}=25\%$. More generally, if $k<N$ there exist exactly $2^{N-k}$ different paths such that their first $k$ vertices lie on a preassigned path $\{ S_{1}, \dots , S_{k} \}$.

Now let us imagine the motion of a particle along a vertical axis as an indicator that shows the cumulative gain at all times of the two players: for $P$, one unit step upward if the coin lands on head, one unit step downward if the coin lands on tail, and the opposite for $Q$. This particle performs a real random walk and its path $\{ S_{1}, \dots , S_{N} \}$ represents the space-time diagram of its motion. In particular, we shall say that

\begin{itemize}
	\item at time $i$ we have a return to the origin if $S_i=0$,
	\item at time $i$ we have a first return to the origin if $\{ S_{1}\neq 0,\dots , S_{i-1}\neq 0, S_{i}=0 \}$,
	\item at time $i$ we have a first passage through $A$ if $\{ S_{1}<A,\dots , S_{i-1}<A, S_{i}=A \}$.
\end{itemize}

Since we will be primarily interested in a possible return to the origin and a return can occur only at even times, as we said, we set $N=2n$. The following proposition lists the probabilities of some possible events we could be interested in.

\begin{proposition}
\label{theorem3}
The following equalities hold:
\begin{equation*}
\begin{split}
a)\quad &  {\cal P}(S_{2n}=0)=\frac{1}{2^{2n}}\binom{2n}{n}\\
b)\quad &  {\cal P}(S_{1}\neq 0,\dots , S_{2n} \neq 0 )=\frac{1}{2^{2n}}\binom{2n}{n}\\
c)\quad &  {\cal P}(S_{1}\geq 0,\dots , S_{2n} \geq 0 )=\frac{1}{2^{2n}}\binom{2n}{n}\\
d)\quad &  {\cal P}(S_{1}\neq 0,\dots , S_{2n-1} \neq 0, S_{2n} = 0 )=\frac{1}{n\cdot 2^{2n-1}}\binom{2n-2}{n-1}\\
e)\quad &  {\cal P}(S_{1}\geq 0,\dots , S_{2n-2} \geq 0, S_{2n-1} < 0 )=\frac{1}{n\cdot 2^{2n-1}}\binom{2n-2}{n-1}\\
\end{split}
\end{equation*}
\end{proposition}

In words: the three events $a)$ a return to the origin takes place at time $2n$, $b)$ no return occurs up to and including time $2n$, and $c)$ the path is non-negative between $0$ and $2n$ have all the same probability $\frac{1}{2^{2n}}\binom{2n}{n}$; whereas, the two events $d)$ the first return to the origin takes place at time $2n$, and $e)$ the first passage through $-1$ occurs at time $2n-1$, have both the same probability $\frac{1}{n\cdot 2^{2n-1}}\binom{2n-2}{n-1}$.

We introduce this notation that will come in very handy in what comes next:
\begin{equation}
\label{eq7}
u_{2n}=\binom{2n}{n}p^n q^n
\end{equation}
which for $p=q=\frac{1}{2}$ becomes
\begin{equation}
\label{eq8}
u_{2n}=\frac{1}{2^{2n}}\binom{2n}{n}.
\end{equation}

In this section, we always refer to the specific case in Eq. (\ref{eq8}). The probabilities of the cases $a)$, $b)$, and $c)$ are then exactly $u_{2n}$ and the probabilities of the cases $d)$, and $e)$ can be expressed as $\frac{1}{2n}u_{2n-2}=u_{2n-2}-u_{2n}=\frac{1}{2n-1}u_{2n}$.

The proof of this proposition is in the appendix. Here we want to use this result to analyze one of the main statements of this reading. To do that, we need this definition: we shall say that the particle that moves according to the random walk spends the time from $k-1$ to $k$ on the positive side if the $k-$th side of its path lies above the $x$-axis, i.e. if at least one of the two vertices $S_{k-1}$ and $S_k$ is positive.

Here it is our main theorem.

\begin{proposition}
Let $p_{2k,2n}$ be the probability that in the time interval from $0$ to $2n$ the particle spends $2k$ time units on the positive side and $2n-2k$ time units on the negative side. Then
\begin{equation}
p_{2k,2n}=u_{2k}\cdot u_{2n-2k}
\end{equation}
\end{proposition} 
\noindent Again the reader can find the complete proof in the appendix. By Eq. (\ref{eq8}), $p_{2k,2n}$ can be given the following equivalent expression
\begin{equation}
p_{2k,2n}=\frac{1}{2^{2k}}\binom{2k}{k}\cdot \frac{1}{2^{2n-2k}}\binom{2n-2k}{n-k}=\frac{1}{2^{2n}}\frac{\binom{n}{k}^2 \binom{2n}{n}}{\binom{2n}{2k}}
\end{equation}

\noindent We want to stress now the importance of this result by means of the following remarks.

\hfill

\begin{remark}
We intuitively feel that the fraction $\frac{k}{n}$ of the total time spent on the positive side is most likely close to $1\over 2$, i.e. that the particle spends half its time above the axis and half its time below it. But will this be true? It may seem paradoxical, but we will see that it is not, and that the opposite is true.

In fact, have a look at the following Table \ref{Table1} which lists the values of the probabilities $p_{2k,2n}$ as a function of $2k$ for $2n=20$ and at the Fig. \ref{fig3} which represents the same probabilities for $2n=20$ and $2n=100$.

\begin{table}[H]
	\centering
	\renewcommand{\arraystretch}{1.2}
\scriptsize
\setlength\tabcolsep{3pt}
\begin{tabular}{||c|c|c|c|c|c|c|c|c|c|c|c||}
\hline
{\bf {$2k$}}                &\ 0 \ &\ 2 \ &\ 4 \ &\ 6 \ &\ 8 \ &\ 10 \ &\ 12 \ &\ 14 \ & 16 \ & 18 \ & 20 \\ \hline
{\bf $p_{2k,2n}$}      & 0.176  & 0.093  & 0.074  & 0.065  & 0.062  &  0.061  & 0.062  & 0.065   & 0.074  & 0.093 & 0.176 \\ \hline
\end{tabular}
\caption{Values of the probabilities $p_{2k,2n}$ that a particle spends a fraction $k/n$ of its time on the positive side, as a function of $2k$ for $2n=20$.}
\label{Table1}
\end{table}

\begin{figure}[H]
\centering
	\subfloat[]{\includegraphics[width=0.45\textwidth]{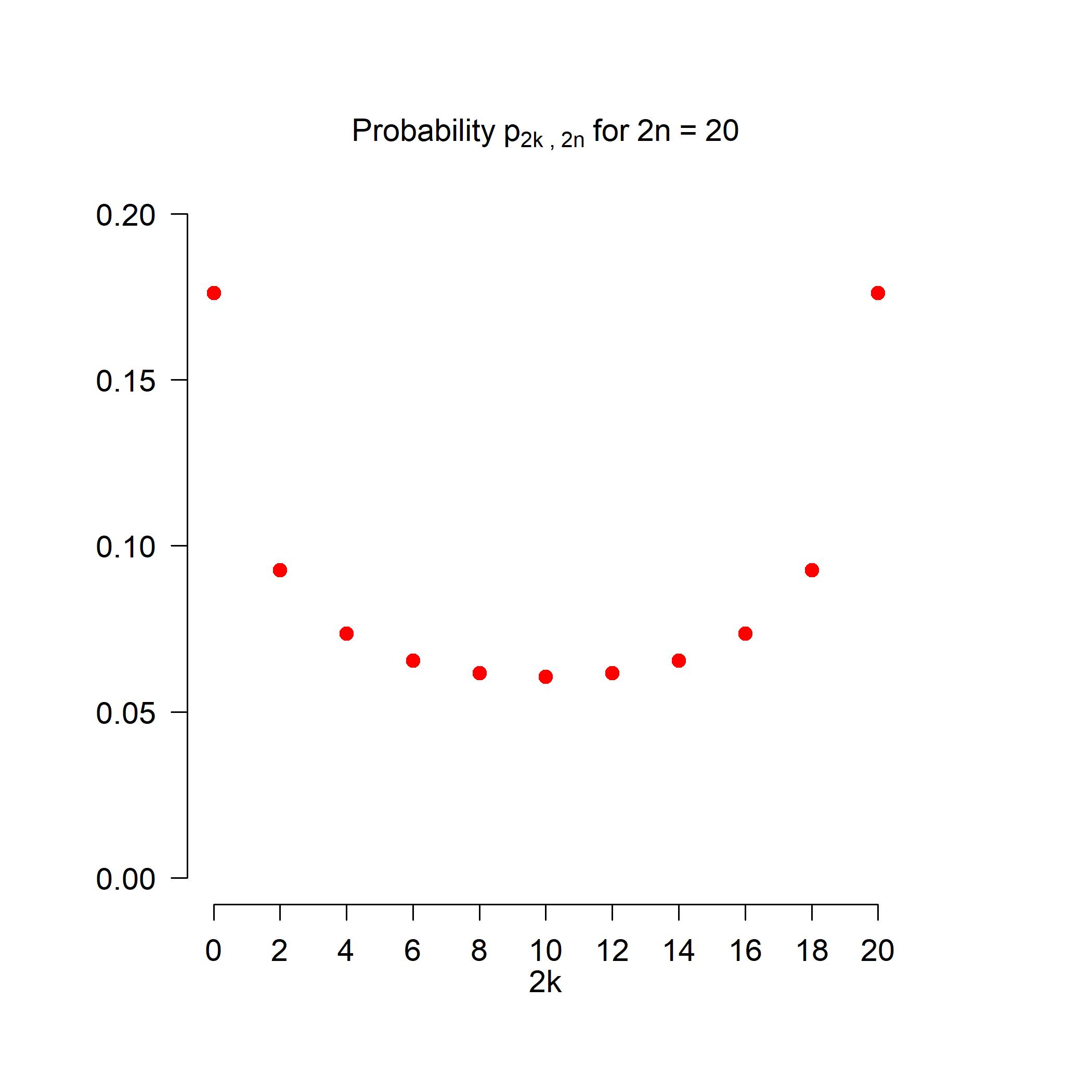}}
	\subfloat[]{\includegraphics[width=0.45\textwidth]{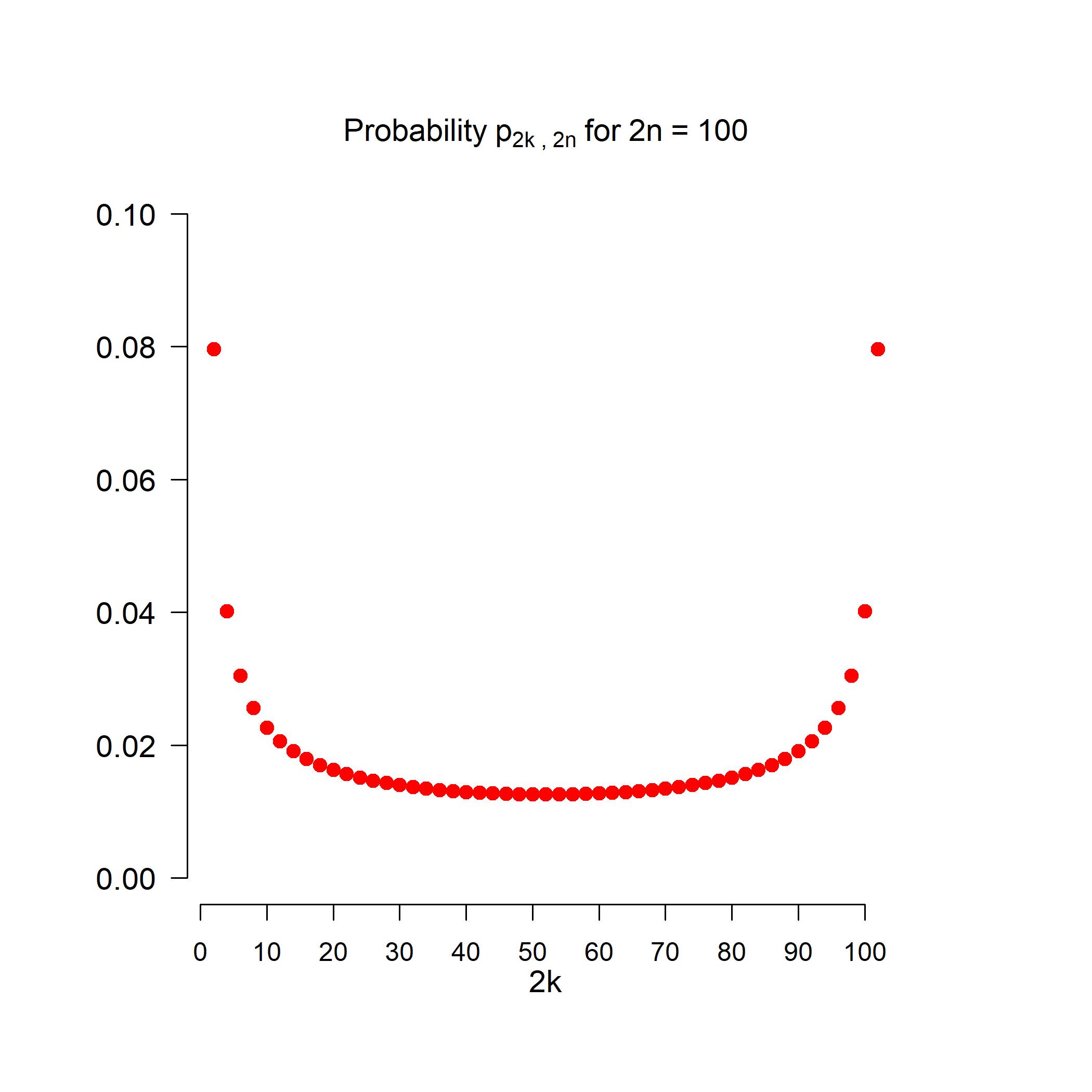}}
	\caption{Plot of the probabilities $p_{2k,2n}$ that a particle spends a fraction $k/n$ of its time on the positive side, as a function of $2k$, (a) for $2n=20$ and (b) for $2n=100$.}
	\label{fig3} 
\end{figure}

As we can see, values close to $\frac{k}{n}={1\over 2}$ are the least probable, whereas the highest probability is assigned to the extreme values $\frac{k}{n}=0$ and $\frac{k}{n}=1$. What does it mean? If we toss a coin $20$ times, the probability $p_{10,20}$ that the particle spends half of its time in the positive side is $0.06=6\%$, whereas the probability that the particle spends all time on the positive side or all time on the negative side is $0.18=18\%$, three times the previous one.
In general, the two extreme probabilities $p_{0,2n}=p_{2n,2n}$ and $p_{n,2n}$ are given by
\begin{equation}
\begin{split}
p_{0,2n} & = p_{2n,2n} = \frac{1}{2^{2n}}\binom{2n}{n} = \frac{1}{2^{2n}} \frac{(2n)!}{(n!)^2}\\
p_{n,2n} & = \frac{1}{2^{2n}}{\binom{n}{n/2}}^2 = \frac{1}{2^{2n}} \frac{(n!)^2}{({n \over 2}!)^4}\\
\end{split}
\end{equation}
and their ratio is given by
\begin{equation}
\frac{p_{0,2n}}{p_{n,2n}}= (2n)! \left( \frac{{n\over 2}!}{n!} \right) ^4
\end{equation}
We can use an important approximation of the factorials for very large values of $n$, called the Stirling approximation, by which we can observe that 
\begin{equation}
\frac{p_{0,2n}}{p_{n,2n}} \sim \frac{\sqrt{\pi n}}{2} \sim \sqrt{n} \qquad {\rm for} \qquad n\to +\infty
\end{equation}
The ratio between the two extreme probabilities increases with the number of trials as $\sqrt n$, which means, for example, that for $2n=200$ trials, the probability of the particle spending all of its time on the positive side is ten times the probability of an equally distributed position above and below the axis. So we begin to realize that if we start winning, we will keep winning, but also that if we start losing, we will keep losing! and this with much greater probability than a more or less balanced alternation of wins and losses! Let us try to restate the same idea in a slightly different way.
\end{remark}

\begin{remark}
The quantity $p_{2k,2n}$ returns the probability that the particle spends exactly $k/n$ of its time on one side, positive or negative. Now let us consider the probabilities that the particle spends \textit{at most} or \textit{at least} $\alpha/n$, for a given $\alpha \in {\mathbb N}$, of its time on one side. The so called \textit{cumulative distribution function} ${\cal P}(2k\leq 2\alpha)$, $\alpha \in {\mathbb N}$, $0\leq \alpha\leq n$ returns this probability and it is defined as
\begin{equation}
\begin{split}
{\cal P}(2k\leq 2\alpha) & =\sum_{k=0}^{\alpha} p_{2k,2n} = \sum_{k=0}^{\alpha} \frac{1}{2^{2n}}\binom{2k}{k}\binom{2n-2k}{n-k} \\
                  & = \frac{1}{2^{2n}}\binom{2n}{n} \sum_{k=0}^{\alpha} \frac{{\binom{n}{k}}^2}{\binom{2n}{k}}= p_{0,2n} \sum_{k=0}^{\alpha} \frac{{\binom{n}{k}}^2}{\binom{2n}{k}} \\
\end{split}
\end{equation}

It represents precisely the probability that the particle spends up to $2\alpha$ time units in the positive side, i.e. {\it at most} $2\alpha$ time units in the positive side and {\it at least} $2n-2\alpha$ time units in the negative side. For instance, for the previous set of $2n=20$ trials the distribution function ${\cal P}(2k\leq 2\alpha)$ is given in the following Table \ref{Table2}

\begin{table}[H]
	\centering
	\renewcommand{\arraystretch}{1.2}
\scriptsize
\setlength\tabcolsep{3pt}
\begin{tabular}{||c|c|c|c|c|c|c|c|c|c|c|c||}
\hline
{\bf {$2\alpha$}}                &\ 0 \ &\ 2 \ &\ 4 \ &\ 6 \ &\ 8 \ &\ 10 \ &\ 12 \ &\ 14 \ & 16 \ & 18 \ & 20 \\ \hline
{\bf {${\cal P}(2k\leq 2\alpha)$}}      & 0.176  & 0.270  & 0.342  & 0.408  & 0.470  &  0.530  & 0.592  & 0.658   & 0.730  & 0.824 & 1.000 \\ \hline
\end{tabular}
\hfill
\caption{Values of the cumulative probabilities ${\cal P}(2k\leq 2\alpha)$ as a function of $2k$ for $2n=20$.}
\label{Table2}
\end{table}
Data in Table \ref{Table2}, and the analogous ones for $2n=100$, are depicted in the Fig. \ref{fig4}.

\begin{figure}[H]
\centering
	\subfloat[]{\includegraphics[width=0.40\textwidth]{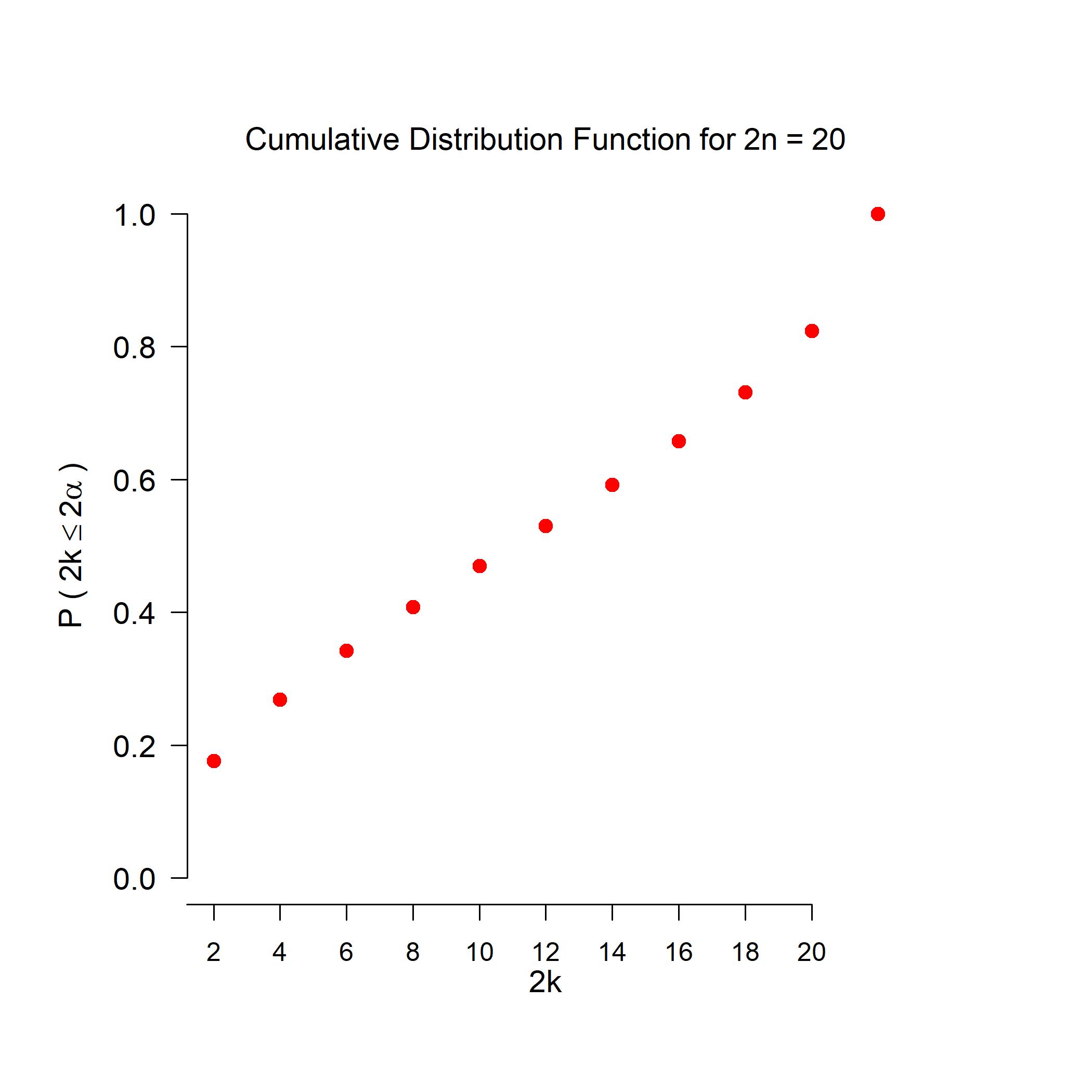}}
	\subfloat[]{\includegraphics[width=0.40\textwidth]{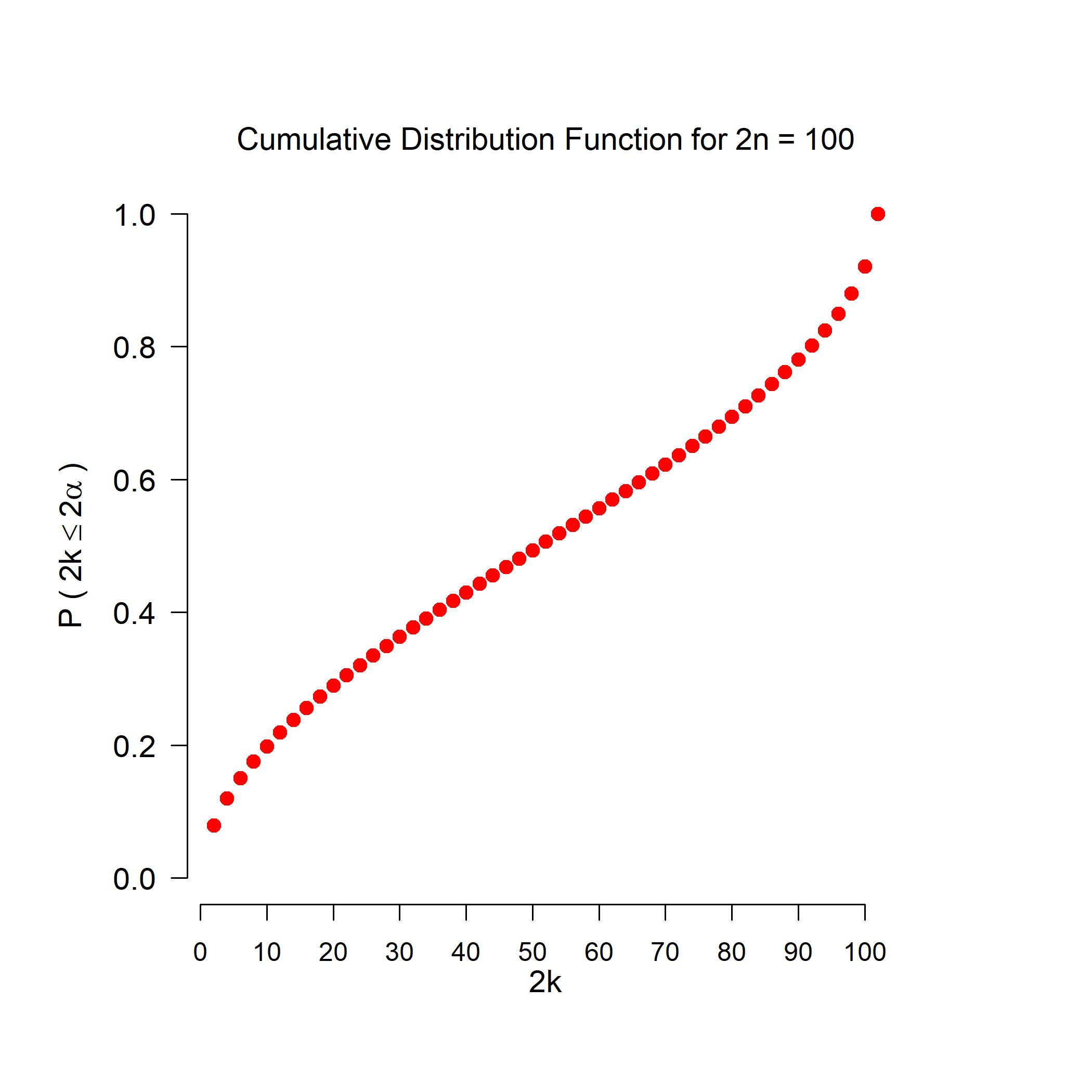}}
	\caption{Cumulative distribution function ${\cal P}(2k\leq 2\alpha)$ as a function of $2k$ (a) for $2n=20$ and (b) for $2n=100$.}
	\label{fig4} 
\end{figure}

The behavior shown in these plots is not coincidental, but represents a very general property of cumulative probability. In fact, these cumulative distribution functions follow a recurrent {\it arcsin} law, as stated by the following proposition.
\begin{proposition}
For $2n\to +\infty$, the cumulative distribution function ${\cal P}(2k\leq 2\alpha)$ can be approximated by
\begin{equation}
{\cal P}(2k\leq 2\alpha)\approx \frac{2}{\pi}{\rm arcsin}\, \sqrt{\frac{\alpha}{n}}
\label{arcsinlaw}
\end{equation}
\end{proposition}

This formula returns with good approximation the probability that the time spent on the positive side is at most, or, equivalently, less than $2\alpha$ with $0\leq 2\alpha\leq 2n$.
The following Fig. \ref{fig5} compares exact results with results expected by the {\it arcsin} model:

\begin{figure}[H]
\centering
	\subfloat[]{\includegraphics[width=0.40\textwidth]{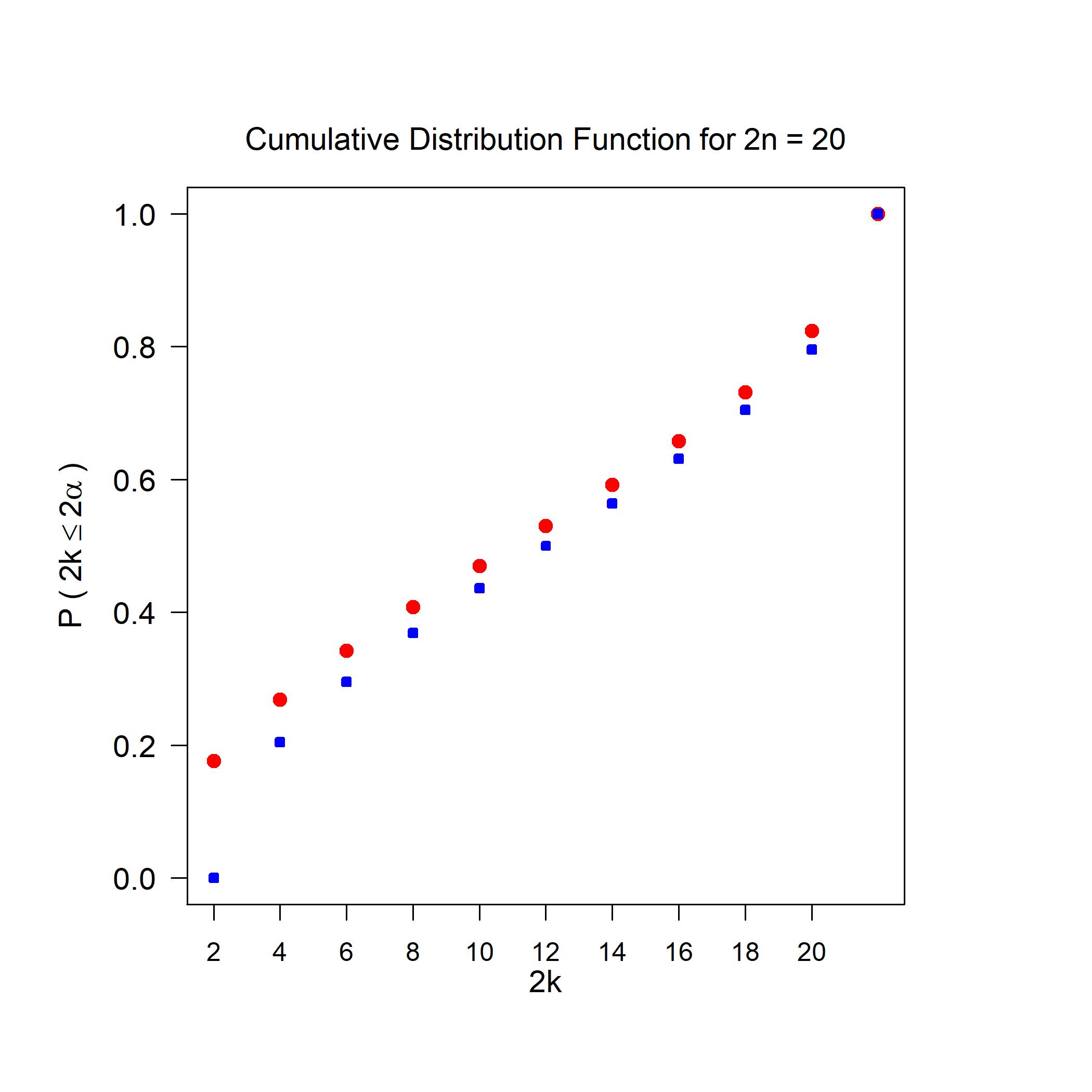}}
	\subfloat[]{\includegraphics[width=0.40\textwidth]{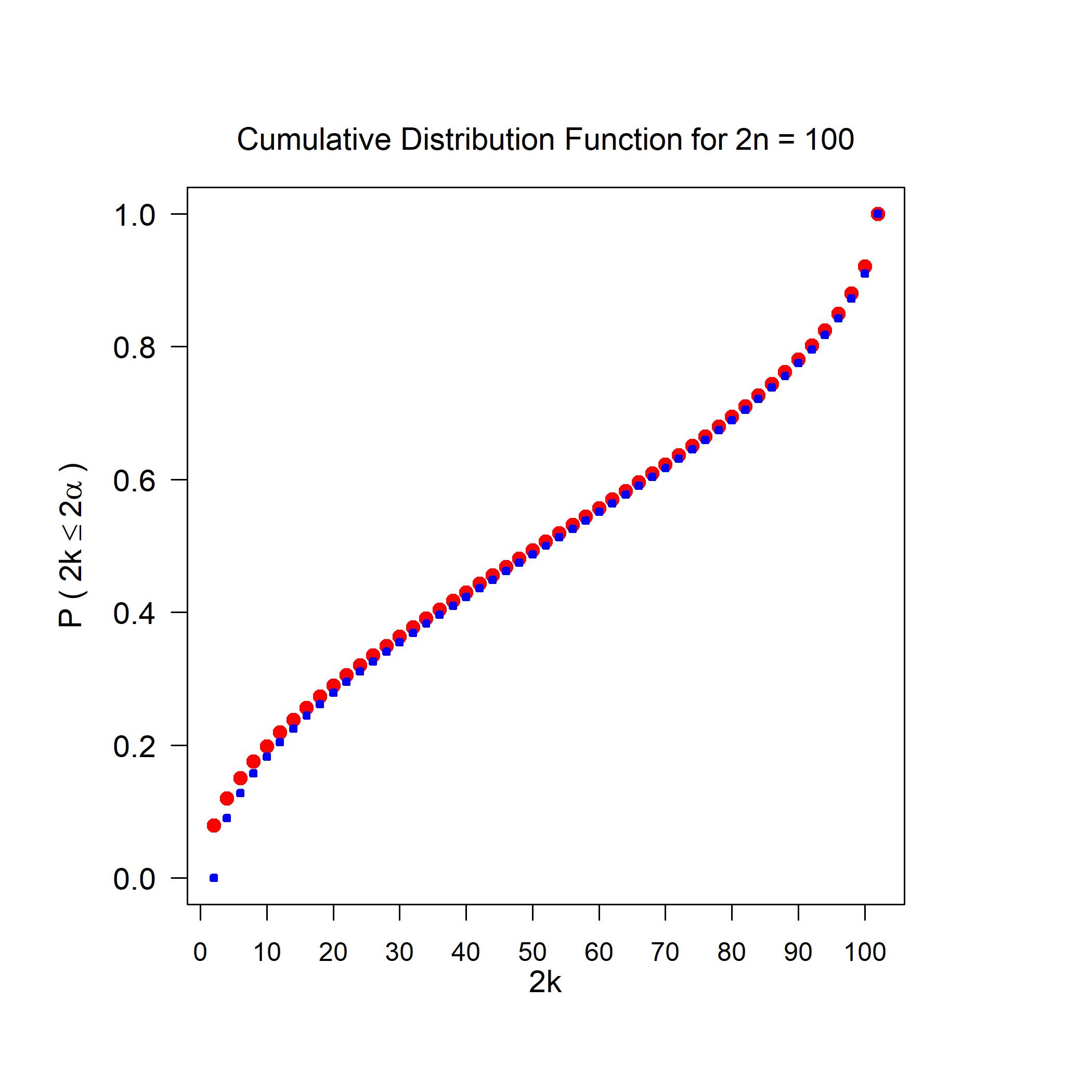}}
	\caption{Comparison between exact data (red circle points) and \textit{arcsin law} (blue square points) (a) for $2n=20$ and (b) for $2n=100$.}
	\label{fig5} 
\end{figure}
Let's be clear about the meaning of formula (\ref{arcsinlaw}). Consider the case $2n=20$ and the value $2\alpha=6$. The cumulative probability is approximately ${\cal P}(2k\leq 6)\approx 40\%$. This is the probability that the particle spends \textit{at most} $\alpha/n=6/20=30\%$ of its time on the positive side and \textit{at least} $1-\alpha/n=14/20=70\%$ of its time on the negative side. In other words, if the particle were a gambler, with a $40\%$ probability it would spend no more than $30\%$ of its time in the game in a winning role. It is a fairly high probability compared to the small fraction of time he plays as a winner! In this case, the gambler-let's say $P$-would be the less fortunate of the two, spending less than half of the game time in a winning role, i.e. he would lose more than he wins. 
\end{remark}

\begin{remark}
Let us consider the event $E_1=$\{The player $P$ is in lead {\it less than} $2\alpha$ days\}. This event can be reformulated as $E_1=$\{$P$ is in lead for $\rho\leq x$ \} where $\rho=\frac{2k}{2n}$ and $x=\frac{2\alpha}{2n}$ are interpreted as fractions of the game time.
This event has probability
\begin{equation}
{\cal P}(E_1)=\frac{2}{\pi}{\rm arcsin}\, \sqrt{x}
\end{equation}
This is the approximate probability that player $P$ spends less time than $x$ in the positive region above the axis, i.e. in the role of winner.
Suppose he is the less fortunate of the two players, so that $x<\frac{1}{2}$. The same can be said of course by swapping the roles of the two players: $E_2=$\{The player $Q$ is in lead for {\it less than} $2\alpha$ days\}$=$\{$Q$ is in lead for $\rho\leq x$ \}. The latter event has the same probability as the former
\begin{equation}
{\cal P}(E_2)=\frac{2}{\pi}{\rm arcsin}\, \sqrt{x}
\end{equation}
We cannot know \textit{a priori} who is the unluckiest player. Consider the event $E={E_1}\cup {E_2}$, that is
\{One of the players is in lead for {\it less than} $2\alpha$ days\}. We are assuming that $x$ is less than half of the playing time so we are precisely considering the point of view of the least fortunate player, whoever he or she may be. The probability of event $E$ is then
\begin{equation}
{\cal P}(E)=2\cdot\frac{2}{\pi}{\rm arcsin}\, \sqrt{x}=\frac{4}{\pi}{\rm arcsin}\, \sqrt{x}
\end{equation}
Suppose $P$ and $Q$ decide to extend the game for a whole year, $2n=365$ days! During this year, we don't know which of the two, but one of the two players is the less fortunate one and is in lead for a fraction $x<\frac{1}{2}$ of the year. Call this probability ${\cal P}(E)=p$. We have
\begin{equation}
p=\frac{4}{\pi}{\rm arcsin}\, \sqrt{x}
\end{equation}
and $x$, that is the {\it maximum} fraction of year in which the less fortunate player is in lead with probability $p$, is given by
\begin{equation}
x=\sin^{2}\left( \frac{\pi}{4}p \right)
\end{equation}
Let's try to give some numbers: if $p=0.05$, then $x=0.00154=0.5625\ {\rm d} =13.50\ {\rm h}$. This means that, with probability $0.05$, i.e. $1$ out of $20$ cases, the less fortunate player will be in lead {\it at most} for $13.50$ hours and the more fortunate player will be in lead for {\it at least} $364$ days and $10.5$ hours. There is a not inconsiderable probability that the unluckiest player will be in the lead for only a tiny fraction of the year! Table \ref{Table3} collects similar results for $x$, expressed in days or hours, for different values of $p$.
\begin{table}[H]
	\centering
	\renewcommand{\arraystretch}{1.2}
\scriptsize
\setlength\tabcolsep{3pt}
\begin{tabular}{||c|c|c|c|c|c|c|c|c|c|c|c|c|c|c||}
\hline
$p$ & \small 0.99 & \small 0.95  & \small 0.90  & \small 0.80  & \small 0.70  & \small 0.60  & \small 0.50  & \small 0.40  & \small 0.30  & \small 0.20  & \small 0.10 & \small 0.05 & \small 0.02 & \small 0.01 \\ \hline
$x$ & \small 179,6  & \small 168,2  & \small 154,0  & \small 126,1 & \small 99,6  & \small  75,2  & \small 53,5  & \small 34,9  & \small 19,9  & \small 8,9 & \small 2,2 & \small 13,5 h & \small 2,2 h & \small 0,5 h \\ \hline
\end{tabular}
\hfill
\caption{Values of the probability $p$ as a function of the fraction of year $x$, expressed in days.}
\label{Table3}
\end{table}
\noindent For example, there is a $50\%$ chance that the least lucky player will be in the lead for a maximum of $53$ days. Let us notice that, when $p\to 0$, $x\to 0$, and, when $p\to 1$, $x\to {1 \over 2}$.
\end{remark}

\begin{remark}
What was shown in the previous remark has an interesting implication. It says that frequently enormously many trials are required before the particle returns to the origin, or that it takes an enormous number of attempts for the less fortunate player to at least break even. This means that the path crosses the $x$-axis very rarely. If we toss a coin for $2n$ times, the number of ties will be proportional to $\sqrt{2n}$. It can be proved that the probability that, within the time $2n$, the particle returns to zero exactly $r$ times is given by
\begin{equation}
{\cal P}_{0}(r,2n)=\frac{1}{2^{2n-r}}\binom{2n-r}{n}
\end{equation}
In particular, $u_{2n}={\cal P}_{0}(0,2n)={\cal P}_{0}(1,2n)>{\cal P}_{0}(2,2n)>{\cal P}_{0}(3,2n)>\dots$. For instance, for $2n=100$, we have ${\cal P}_{0}(0,100)={\cal P}_{0}(1,100)=7.96\%$, ${\cal P}_{0}(5,100)=7,17\%$, ${\cal P}_{0}(10,100)=4.84\%$, ${\cal P}_{0}(20,100)=0.73\%$, and ${\cal P}_{0}(30,100)=0.014\%$.

We have hitherto assumed that the two probabilities $p$ and $q$ were equal. This is indeed the case in a coin toss but we can imagine a random walk produced by a process in which the two choices are not equiprobable. We now want to extend our argument to the case in which $p\neq q$ but, in order to do that, we need some further results about the so called {\it ruin problem}, which is the topic of the next section.
\end{remark}

\section{The Ruin Problem}
\label{section4}

In the fair game of the previous section, $p$ and $q$ were interpreted as the probabilities of getting heads or tails on a coin flip. All in all, this idea is fairly intuitive because the two events are reasonably equiprobable. In the ballot problem $p$ and $q$ were numbers, the number of actual outcomes of a ballot on $N$ repeated trials. However, the ratios $p/N$ and $q/N$ are frequencies, and, if $N$ is large enough, they are the probabilities of the two different outcomes. From now on, $p$ and $q$ will directly represent the two probabilities that one of two possible results will come out at each repeated trial, that is, that a given random variable will symbolically take on the value $+1$ or $-1$. If each draw or roll is independent of the previous one, then in no way does one result affect the next. Technically, it is said that we have a set of $N$ \textit{independent and identically distributed random variables}.

We are now in a position to give a more general and more formal definition of the classical random walk. The random walk in one dimension is the stochastic process described by the motion of a particle whose position, at each step, is
\begin{equation}
S_n=S_0+\sum_{k=1}^{n}X_k
\end{equation}
where $X_k$ are independent and identically distributed (i.i.d.) random variables such that ${\cal P}(X_k=1)=p$ and ${\cal P}(X_k=-1)=q$, with $p+q=1$. We can assume $S_0=0$, without losing anything in our reasoning.

What we're going to study now is the behavior of the random walk in the presence of two barriers, that is, two different values of $S_n$, let's say a positive $y=A$ and a negative $y=-B$, so that the random walk stops when it reaches one of these two values. For example, suppose one of two players $P$ and $Q$ stops gambling when he or she gets $A$ wins or $B$ losses in two-outcome game. This is called \textit{ruin problem} and we want to find out the probability to reach $y=A$ before $y=-B$ or viceversa and the mean duration of this process. We will start from the case $p=q=\frac{1}{2}$ which is called {\it unbiased} since there is no drift toward the positive values or the negative values. Then, we will move to the {\it biased}, and more general, case, in which the two probabilities are different, $p\neq q$. In both cases, we will study similar propositions but their consequences will be very different.

\begin{figure}[H]
	\centering
	\includegraphics[width=0.8\textwidth]{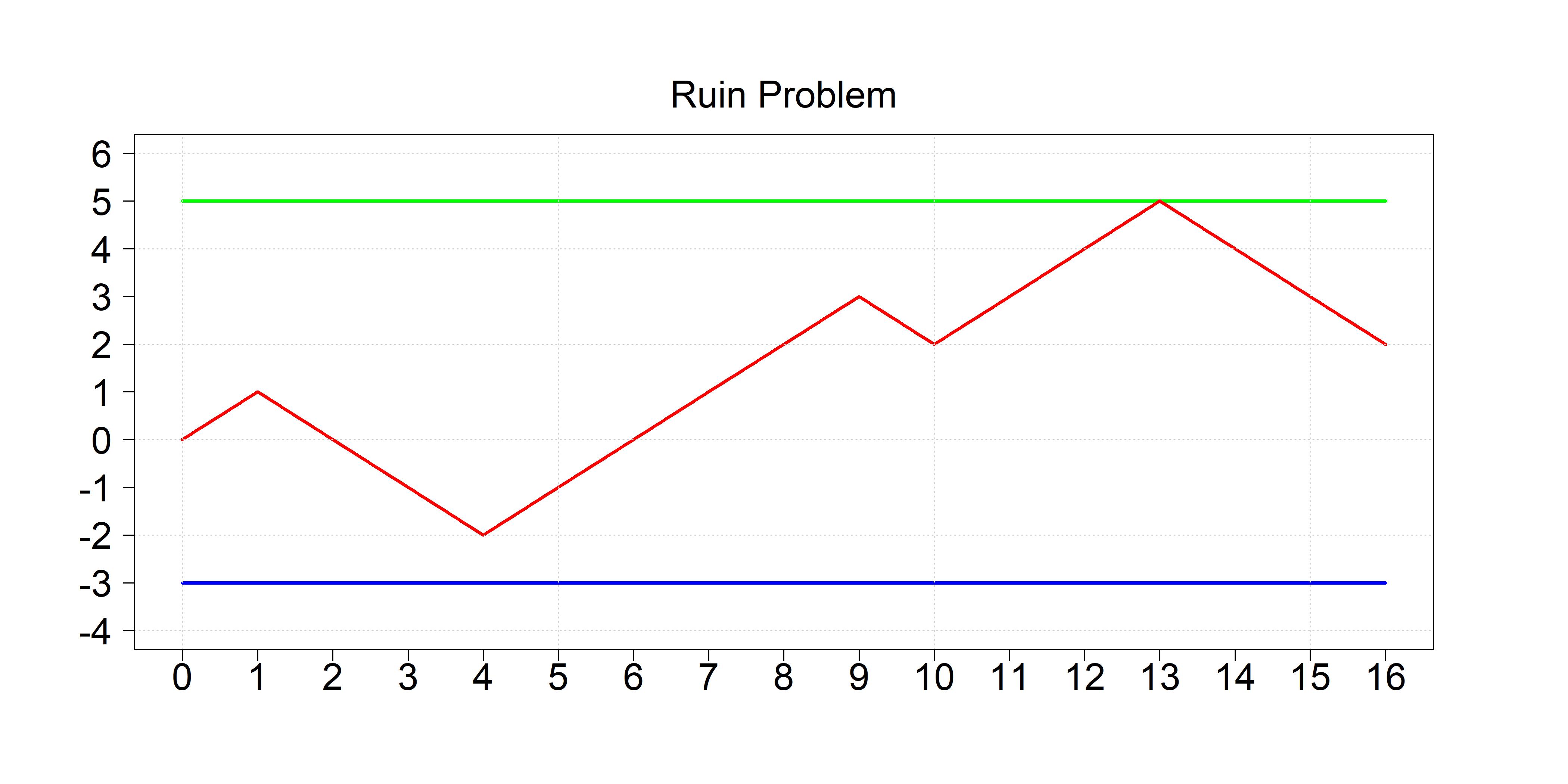}
	\caption{An illustrative example of the ruin problem for an unbiased random walk, with $A=5$ and $B=3$.}
	\label{fig6}
\end{figure}

{\bf The unbiased random walk.} Let us immediately state the following proposition.
\begin{proposition}
Let $A, B \in {\mathbb R}$, $A, B >0$ and $\tau=\min \{ n\geq 0: S_n =A\ {\rm or}\ S_n=-B \}$  be the first passage time through $A$ or $-B$. Then
\begin{equation}
{\cal P}(S_\tau =A) =\frac{B}{A+B}, \quad {\cal P}(S_\tau =-B) =\frac{A}{A+B} \quad {\rm and} \quad {\mathbb E}[\tau]=AB
\end{equation}
\label{theorem6}
\end{proposition}
The quantity ${\cal P}(S_\tau =A)$ represents the probability that the random walk reaches $A$ before $-B$ and similarly for ${\cal P}(S_\tau =-B)$; ${\mathbb E}[\tau]$ is the expected mean value of the time or of the number of steps needed to reach one of the two barriers. Let's have a look at the Fig. \ref{fig6}, where $A=5$ and $B=3$.

According to Proposition \ref{theorem6}, the upper barrier in $A$ is reached before the lower one in $B$ with probability ${\cal P}(S_\tau =A)=\frac{3}{8}= 37.5\%$; similarly, ${\cal P}(S_\tau =-B)=\frac{5}{8}=62.5\%$. The expected mean time is ${\mathbb E}[\tau]=10$ units time or steps.

This is a well-known result and it can be easily extended to the case in which the particle does not start from $(0,0)$ but from a point $(0,k)$, with $k\in[-B,A]$. For instance, the expected duration becomes ${\mathbb E}[\tau]=(A-k)(B+k)$.

{\bf The biased random walk.} Now, let us change the game. Suppose we roll a die and we take one step forward if we get $5$ or $6$ and one step back if we get $1$ or $2$ or $3$ or $4$, so that $p=\frac{1}{3}$ and $q=\frac{2}{3}$. 
In technical language we should say that we are extending the previous argument to the case of a sequence of Bernoulli trials with probability $p$ of success and a probability $q$ of failure, with in general $p\neq q$, i.e. $p,q\neq {1\over 2}$, $p+q=1$. In other terms, we have $X_k$ i.i.d. random variables such that
\begin{equation}
{\cal P}(X_k=1)=p\qquad {\cal P}(X_k=-1)=q
\end{equation}
and we want to consider the \textit{biased} random walk defined by
\begin{equation}
S_n=S_0+\sum_{k=1}^{n}X_k
\end{equation}
Let us assume again $S_0=0$. In general we call \textit{unbiased} the random walk if $p=q$ and \textit{symmetric} if $A=B$. Here now is the extension of Proposition \ref{theorem6} to the biased case.

\begin{proposition}
Let $S_n$ be a biased classical random walk, $A,\, B \in {\mathbb R}$, $A,\, B\geq 0$ and $\tau=\min \{ n\geq 0: S_n =A\ {\rm or}\ S_n=-B \}$ be the first passage time through $A$ or $-B$. Then we have
\begin{equation}
\label{eq26}
{\cal P}(S_{\tau}=A)=\frac{1-\big(\frac{q}{p}\big)^{B}}{1-\big(\frac{q}{p}\big)^{A+B}}, \quad
{\cal P}(S_{\tau}=-B)=\frac{1-\big(\frac{p}{q}\big)^{A}}{1-\big(\frac{p}{q}\big)^{A+B}}
\end{equation}
and
\begin{equation}
\begin{split}
{\mathbb E}[\tau]&= \frac{B}{q-p}-\frac{A+B}{q-p}\cdot \frac{1-\big(\frac{q}{p}\big)^{B}}{1-\big(\frac{q}{p}\big)^{A+B}}\\ &= \frac{1}{q-p}\left[ A\frac{1-\big(\frac{q}{p}\big)^{B}}{1-\big(\frac{q}{p}\big)^{A+B}}-B \frac{1-\big(\frac{p}{q}\big)^{A}}{1-\big(\frac{p}{q}\big)^{A+B}} \right]
\end{split}
\end{equation}
\label{theorem7}
\end{proposition}
The proof of propositions \ref{theorem6} and \ref{theorem7} is not elementary, and the reader who is not interested in rather technical details may omit their study. However, we want to bring attention to a substantial set of consequences, which we gather in the following remarks.
\begin{remark}
\label{remark5}
First, the unbiased random walk is a limiting case of the biased one for $p,q\to \frac{1}{2}$. This is proved in the appendix.
\end{remark}

\begin{remark}
If $q=0$ and $p=1$, we have ${\cal P}(S_\tau=A)=1$, ${\cal P}(S_\tau=B)\sim \left( \frac{q}{p}\right) ^B \rightarrow 0$ and ${\mathbb E}[\tau]=(A+B)\cdot 1 -B=A$. This means that we get $A$ in $A$ steps. Similarly, if $q=1$ and $p=0$, we have ${\cal P}(S_\tau=A)\sim \left( \frac{p}{q}\right) ^A \rightarrow 0$, ${\cal P}(S_\tau=B)=1$ and ${\mathbb E}[\tau]\sim (-1)\left[A\left(\frac{p}{q}\right)^A-B\right]\rightarrow B$. That is, we get $-B$ in $B$ steps.
\end{remark}

\begin{remark}
Let's go back to the ${\cal P}(S_\tau=A)$ and ${\cal P}(S_\tau=-B)$. They are given by Eqs. \ref{eq26}, and, written as functions of $\rho= q/p$, they become
\begin{equation}
{\cal P}(S_{\tau}=A)=\frac{1-\rho^{B}}{1-\rho^{A+B}}, \quad {\cal P}(S_{\tau}=-B)=\frac{1-\rho^{-A}}{1-\rho^{-(A+B)}}
\end{equation}
In the symmetric case, i.e. for $A=B$, they are simply
\begin{equation}
{\cal P}(S_{\tau}=A)=\frac{1}{1+\rho^{A}}, \quad {\cal P}(S_{\tau}=-A)=\frac{1}{1+\rho^{-A}}
\end{equation}
It is easy to show that they are equal only for $\rho=1$; more in general, their behavior as functions of $\rho$ is shown in the Fig. \ref{fig7} for two different values of $A=B$.

\begin{figure}[H]
\centering
	\subfloat[]{\includegraphics[width=0.45\textwidth]{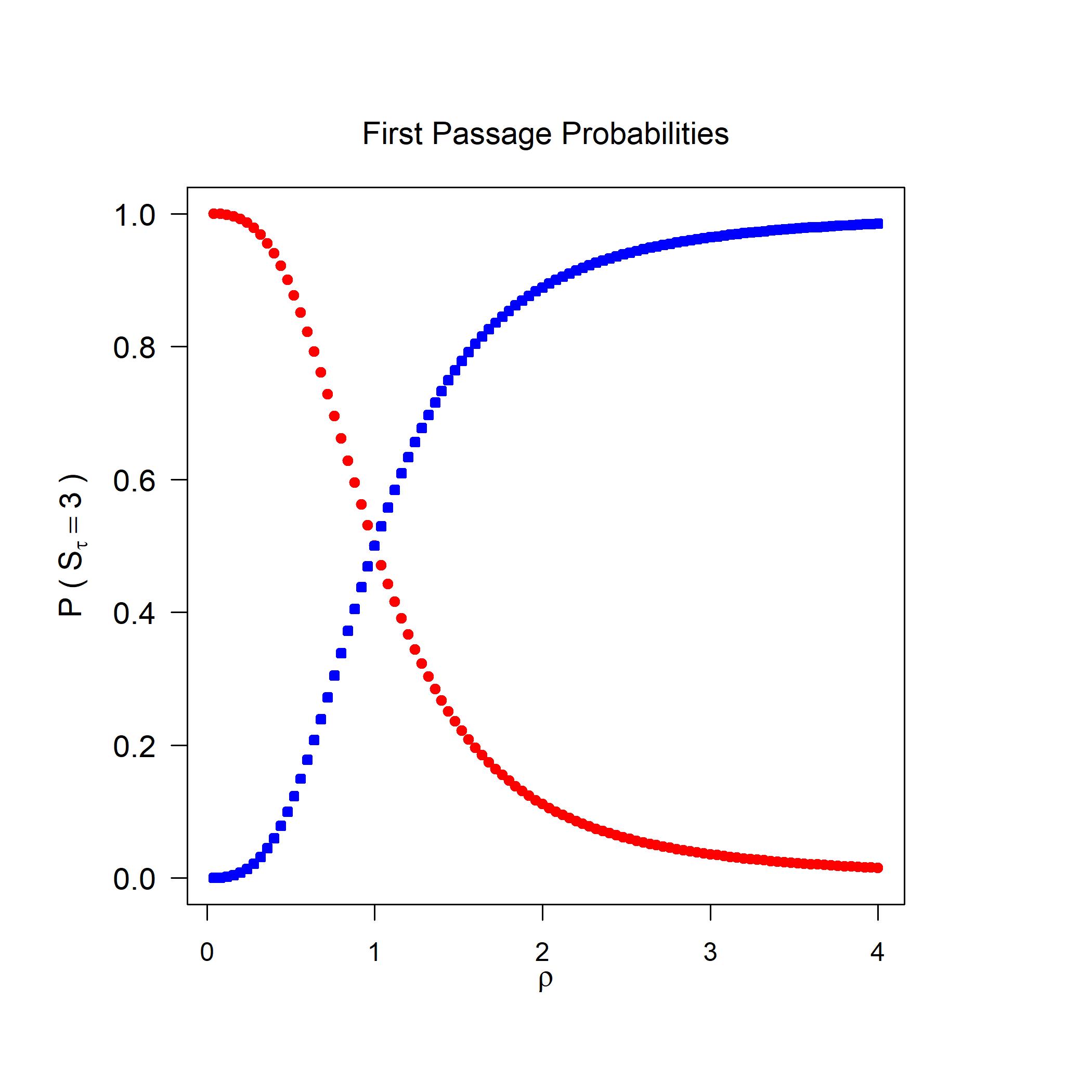}}
	\subfloat[]{\includegraphics[width=0.45\textwidth]{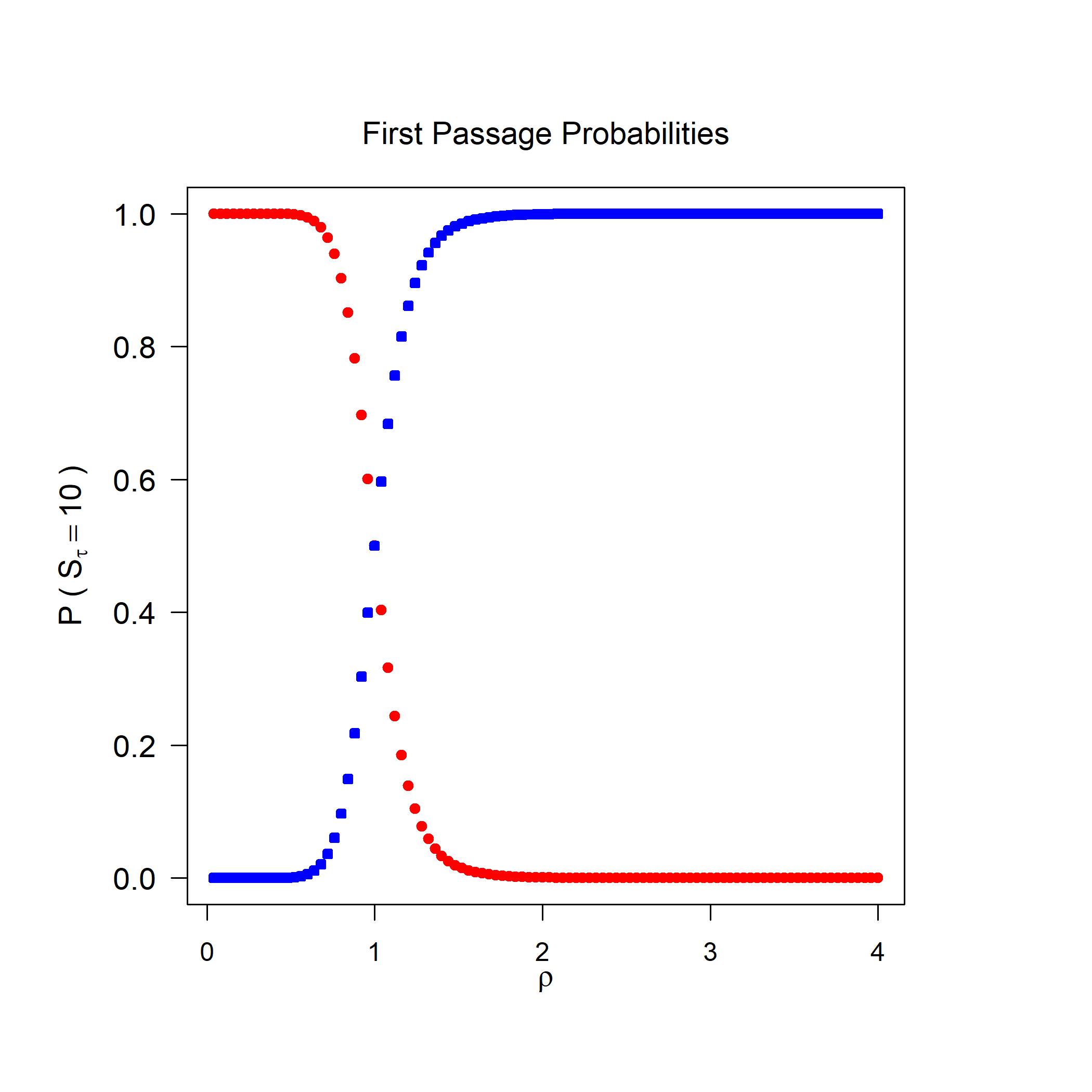}}
	\caption{Plot of the function ${\cal P}(S_{\tau}=A)=\frac{1}{1+\rho^{A}}$ in red circle points and of the function ${\cal P}(S_{\tau}=-A)=\frac{1}{1+\rho^{-A}}$ in blue square points (a) for $A=B=3$ and (b) for $A=B=10$.}
	\label{fig7} 
\end{figure}

The farthest are the barriers with respect to the origin, the quickest is the change in probabilities. For instance, for $A=B=3$ we need roughly speaking a ratio $\rho=4$ to be almost sure not to get the gain at $+3$, whereas for $A=B=10$ it is enough a ratio $\rho=1.5$ to be almost sure not to reach our success at $+10$. Let us stress this point. A high value of $\rho$ means that the probability $q$ of a move down predominates. So when $\rho$ is high, it becomes increasingly difficult to win. In fact, if $\rho$ is about $4$ in the case $A=B=3$, it becomes practically impossible to win. In the case $A=B=10$ a lower ratio $\rho$ between $q$ and $p$ is enough to make it practically impossible to win. In fact, it is reduced to $\rho=1.5$.

If our strategy is to stop when we win $10$ dollars or when we lose $10$ dollars, then the probability of reaching the win first before losing our budget becomes very very small as soon as the game becomes unfavorable to us. For example, it only takes a $\rho=2$ to be virtually certain of losing \$10 before winning \$10!
The following Table \ref{Table4} can give us a deeper insight in the numbers involved.
\begin{table}[H]
	\centering
	\renewcommand{\arraystretch}{1.0}
\scriptsize
\begin{center}
\scriptsize
\begin{tabular}{||c||c|c||c|c||}
\hline \hline
\rule[-4mm]{0mm}{10mm}                          & \multicolumn{2}{c||}{$A=3$} &\multicolumn{2}{c||}{$A=10$}   \\ \hline
\rule[-4mm]{0mm}{1cm}                   		   & ${\cal P}(S_{\tau}=3)$ & ${\cal P}(S_{\tau}=-3)$ & ${\cal P}(S_{\tau}=10)$ & ${\cal P}(S_{\tau}=-10)$  \\ \hline
\rule[-4mm]{0mm}{1cm} $\rho={51 \over 49}$     & $47\%$  & $53\%$  & $40\%$  & $60\%$ \\ \hline
\rule[-4mm]{0mm}{1cm} $\rho={55 \over 45}$     & $35\%$  & $65\%$  & $11\%$  & $89\%$ \\ \hline
\rule[-4mm]{0mm}{1cm} $\rho={60 \over 40}$     & $23\%$  & $77\%$  & $2\%$   & $98\%$ \\ \hline \hline
\end{tabular}
\end{center}
\hfill
\caption{Values of the functions ${\cal P}(S_{\tau}=A)=\frac{1}{1-\rho^{A}}$ and ${\cal P}(S_{\tau}=-A)=\frac{1}{1+\rho^{-A}}$ for $A=3$ and $A=10$.}
\label{Table4}
\end{table}
As we can see from Table \ref{Table4}, the probability of win for a gambler in an unfair game depends on the value of $A$. Suppose that $A$ also represents the final payout if it is reached.
If $\rho={51 \over 49}$, i.e. if he has only the $1\%$ more to loose than to win, the probability to win $3$ dollars is $47\%$, whereas the probability to loose $3$ dollars is $53\%$ but the probability to win $10$ dollars decreases to $40\%$, whereas the probability to loose $10$ dollars increases to $60\%$. This is much more conspicuous if we reduce the probability of win of the $10\%$, that is for $\rho={60 \over 40}$. In this case, the probability to win $3$ dollars is $23\%$, whereas the probability to loose $3$ dollars is $77\%$ and the probability to win $10$ dollars decreases to $2\%$, whereas the probability to loose $10$ dollars increases to $98\%$. An imbalance of only $10\%$ in favor of the casino counter produces a collapse in the value of the odds of winning to $2\%$!

\begin{figure}[H]
	\centering
	\includegraphics[width=0.7\textwidth]{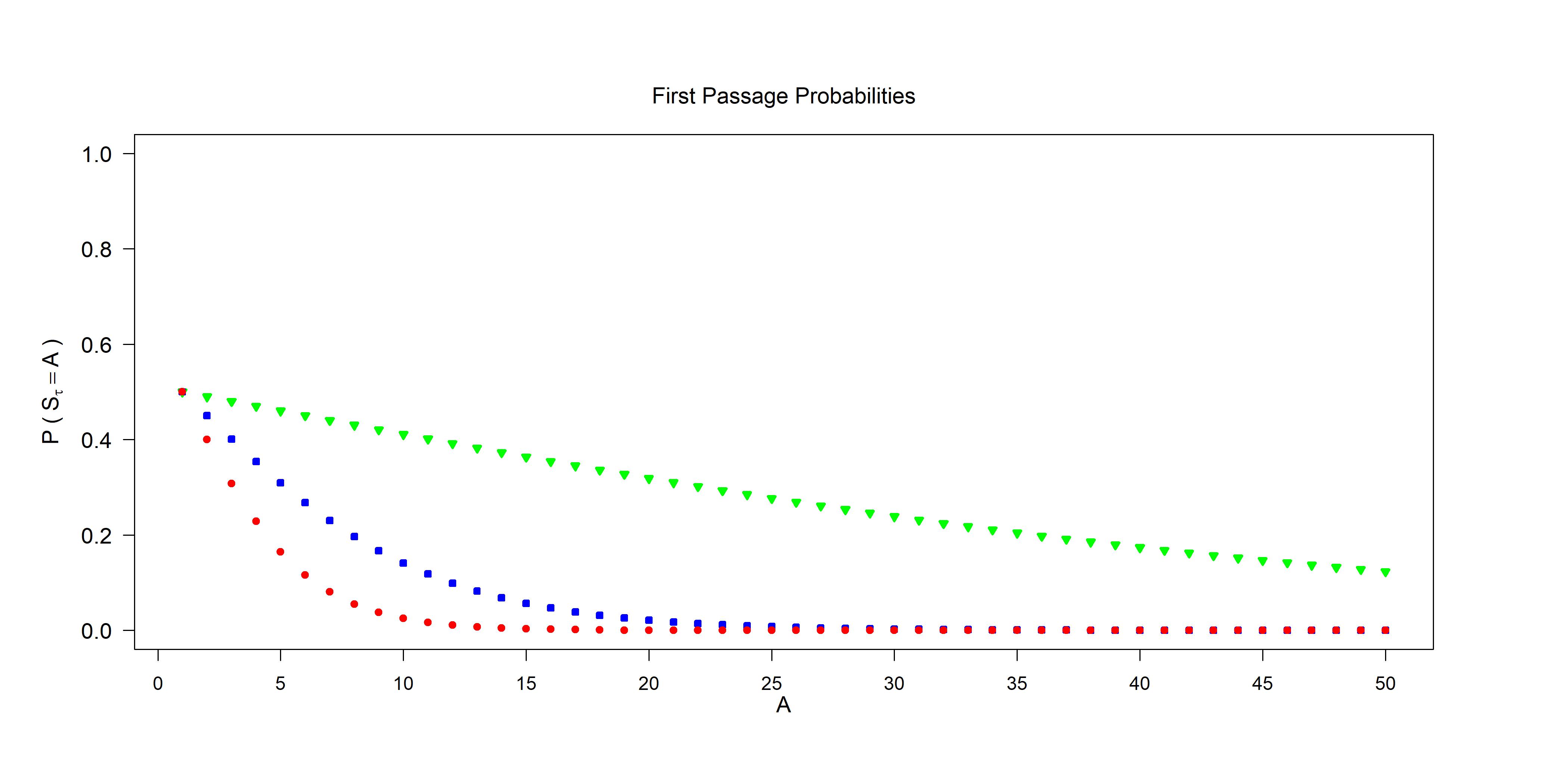}
	\caption{Values of ${\cal P}(S_{\tau}=A)=\frac{1}{1-\rho^{A}}$ as a function of $A$ for three different values of $\rho$: green triangle points $\rho=51/49$; (b) blue square points $\rho=55/45$; (c) red circle points $\rho=40/60$.}
	\label{fig8}
\end{figure}
\end{remark}

It must always be kept in mind that in an unfair game, the probability of losing a given amount of money is much higher than the probability of winning the same amount of money, and the larger the amount, the more noticeable the difference between these probabilities. Fig. \ref{fig8} shows clearly how the probability of win decreases with $A$ for different values of $\rho$.

\begin{remark}
Let us focus now on the expected duration ${\mathbb E}[\tau]$ and study its behavior as a function of $\rho=q/p$ and of the barrier values $A$ and $B$. 
When $A=B$, i.e. with a symmetric barrier, the expected duration becomes
\begin{equation}
	{\mathbb E}[\tau]=A\, \frac{1+\rho}{1-\rho}\cdot \frac{1-\rho^A}{1+\rho^A}
\end{equation}
and the Fig. \ref{fig9} shows its shape as a function of $\rho$ for three different values of $A$:
\begin{figure}[H]
	\centering
	\includegraphics[width=0.9\textwidth]{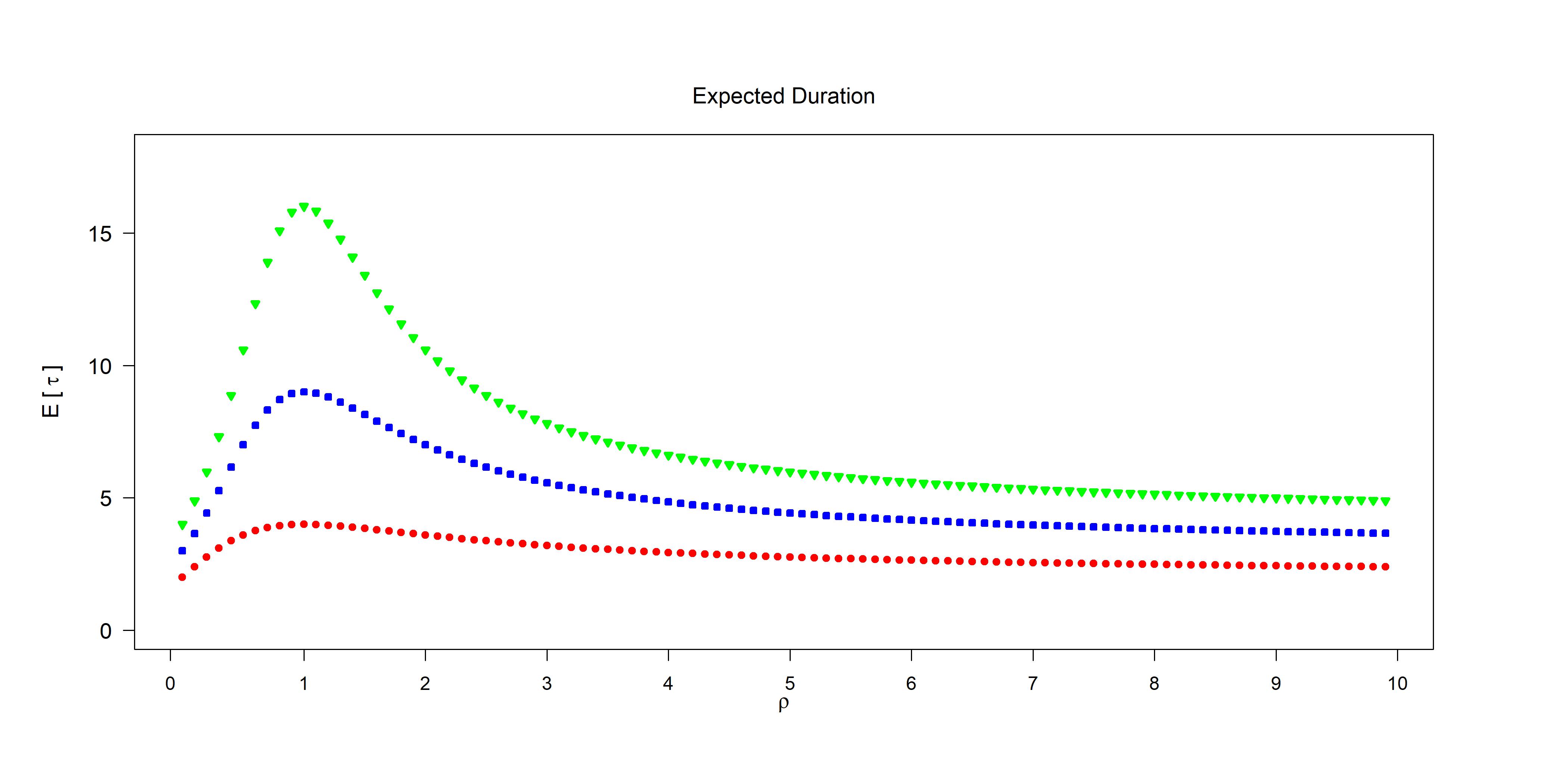}
	\caption{${\mathbb E}[\tau]$ as a function of $\rho$ for three different values of $A$ in the symmetric case: green triangle points $A=4$, blue square points $A=3$, and red circle points $A=2$.}
	\label{fig9}
\end{figure}
It can be observed from the figure--and can also be proved--that each plot has a maximum for $\rho=1$, i.e. $p=q$, and that its value, as expected, is $A^2$. Of course, the plot has an intercept at the origin at $A$ and goes asymptotically to $A$ as $\rho\to+\infty$. In the asymmetric case, the expected duration takes the form
\begin{equation}
{\mathbb E}[\tau]= \frac{\rho+1}{\rho-1}\left[ A\frac{1-\rho^{B}}{1-\rho^{A+B}}-B \frac{1-\rho^{-A}}{1-\rho^{-(A+B)}} \right]
\end{equation} 
 
The behavior is not very dissimilar to the previous one. The intercept at the origin still coincides with the value of $A$ and the asymptotic value with the value of $B$. The curve attains a maximum at a value of $\rho>1$ if $A<B$ and $\rho<1$ if $A>B$. The value of the maximum expected duration grows with $A$ and with $B$. Let us finally observe that the value of the expected duration is perfectly symmetric with respect to an exchange of the values of $A$ and $B$.
\end{remark}

\section{Back again to the origin}
\label{section5}
In the previous discussion, we gained more insight into the biased case, where the probabilities of winning and losing are generally different.
What we want to do now is to go back to the problem of returning to the origin for a random walk, analyzed in the section \ref{section2}, in the case ${\cal P}(X_k=1)=p \neq q= {\cal P}(X_k=-1)$. Studying the return to the origin of our random particle equals focusing on the event $E=$\{The cumulative number of successes and failures are equal\}. Of course, if at the $k$-th trial the cumulative numbers of successes and failures are equal, then $k$ must be an even number, $k=2n$ and $n$ trials must have resulted in success, the other $n$ in failure.

The event $E$ occurs at the $2n$-th trial when $S_{2n}=0$. Moreover, as before, the event $E=$\{The first return to the origin occurs at the $2n$-th trial\} is defined by the aggregate of sequences such that $S_{1}\neq 0, S_{2}\neq 0,\dots , S_{2n-1} \neq 0, S_{2n} = 0$. In this case, the probability of the first few terms can be easily found by direct computation:
\begin{equation*}
\begin{split}
& {\cal P}(S_{1}\neq 0, S_{2}= 0)=2pq  \\
& {\cal P}(S_{1}\neq 0, S_{2}\neq 0, S_{3} \neq 0, S_{4} = 0)=2p^2q^2 \\
& {\cal P}(S_{1}\neq 0, S_{2}\neq 0,\dots , S_{5} \neq 0, S_{6} = 0)= 4p^3q^3 \\
& {\cal P}(S_{1}\neq 0, S_{2}\neq 0,\dots , S_{7} \neq 0, S_{8} = 0)= 10p^4q^4  \\
& {\cal P}(S_{1}\neq 0, S_{2}\neq 0,\dots , S_{9} \neq 0, S_{10} = 0)= 28p^5q^5  \\
\end{split}
\end{equation*}
and so on. For example, if $S_{2}= 0$, the first step could take us up with probability $p$ and down with probability $q$, but if we want to get back to $0$ immediately, we need a second step down with probability $q$ or up with probability $p$, respectively. Summing up the probabilities of two paths we get $2pq$. It is a useful exercise to attempt to obtain the other probabilities listed above. Here we want to generalize Proposition \ref{theorem3}, that we studied in the unbiased case, to the case $p\neq q$.
\begin{proposition}
The following equalities hold:
\begin{equation*}
\begin{split}
a)\quad &  {\cal P}(S_{2n}=0)=\binom{2n}{n}p^n q^n\\
b)\quad &  {\cal P}(S_{1}\neq 0,\dots , S_{2n-1} \neq 0, S_{2n} = 0 )=\frac{2}{n}\binom{2n-2}{n-1}p^n q^n\\
\end{split}
\end{equation*}
\end{proposition} 
In other words, by Eq. (\ref{eq8}) $u_{2n}=\binom{2n}{n}p^n q^n$, the two events $a)$ a return to the origin takes place at time $2n$, and $b)$ the first return to the origin takes place at time $2n$, have probabilities given by 
\begin{equation*}
\begin{split}
a)\quad &  {\cal P}(S_{2n}=0)=u_{2n}\\
b)\quad &  {\cal P}(S_{1}\neq 0,\dots , S_{2n-1} \neq 0, S_{2n} = 0 )=\frac{2pq}{n}\, u_{2n-2}\\ 
& =4pq\, u_{2n-2}-u_{2n}=\frac{1}{2n-1}u_{2n}\\
\end{split}
\end{equation*}

For instance, when $p=\frac{1}{3}$ and $q=\frac{2}{3}$ and we take $2n=20$, the values of the two probabilities are ${\cal P}(S_{20}=0)=5,43\%$ and ${\cal P}(S_{1}\neq 0,\dots , S_{19} \neq 0, S_{20} = 0 )=0,29\%$. Let us notice that the last probability is $1/19$ of the first one.

{\bf A matter of transient or persistent states.} We now want to focus on a subtle issue. The event $E$ related to a return of the random walk to the origin in any of the $2n$-trials is interpreted as the fact the we loose all what we gained or viceversa we re-gain all what we lose. This event could happen a number of times after we started the game.
Let us try to formulate the issue in simple words. How many times do we expect the event $E$ to occur if the number of repeated trials increases? If we imagine repeating the trials \textit{ad infinitum}, what can we say about the expected number of times $E$ will occur?
Although it may seem strange, we want to show that the event $E$ will occur a finite number of times even if we repeat the trials \textit{ad infinitum}. The return to the origin cannot happen infinitely many times.

We now state a theorem. This theorem could appear difficult to understand at first reading. What is more important to us is the subsequent argument, which is more of a proof of the theorem and, for this reason, we do not postpone to the appendix.

Let us denote $E_{2n}$=\{A return to origin occurs at time $2n$\}$=\{S_{2n}=0\}$. Then

\begin{proposition}
The probability that infinitely many events $E_{2n}$ occur is $0$, or, equivalently, ${\cal P}( \limsup_{n\to\infty} E_{2n})=0$. 
\end{proposition}

We will explain the meaning of this proposition gradually.

In Eq. (\ref{eq8}) we defined $u_{2n}$ which is exactly the probability that a return to the origin occurs at time $2n$. Here we are interested in a great number of repeated trials, so we suppose $2n$ big enough. Let us start by using Stirling's formula $\binom{2n}{n}\sim \frac{4^n}{\sqrt{\pi n}}$, for $n\to +\infty$, in order to write
\begin{equation}
u_{2n}\sim \frac{(4pq)^{n}}{\sqrt{\pi n}}
\end{equation}
Let us notice that\footnote{This is related to the fact that the product of two numbers whose sum is constant is maximum if the two numbers are equal.}, if $p\neq q$, then $4pq<1$ so that $u_{2n}\to 0$ and $\sum_{n} u_{2n}$ converges faster than the geometric series with ratio $4pq$. When $p=q=\frac{1}{2}$, $u_{2n}\to 0$ but $\sum u_{2n}$ diverges as it is asymptotic to a divergent generalized harmonic series.\footnote{Let us observe that $\{ u_{2n}\}$ is not a probability distribution since $\sum u_{2n}$ can be greater than $1$. For instance, in our previous example, we have $\sum u_{2n}=3$. Instead, we can interpret $u_{2n}$ as the expectation of a random variable which equals $1$ or $0$ according to weather a return to the origin does or does not occur at $2n$-trial. Hence $\sum_{j=0}^{2n} u_{j}$ is the expected number of occurrences of such a return in $2n$ trials.} In the first case,  when it converges, we can compute the value of the sum of the series. In fact, by setting $x=-4pq$ and $\alpha=-\frac{1}{2}$, and by using the following combinatorial identity
\begin{equation}
\binom{2n}{n}=\binom{-\frac{1}{2}}{n}\cdot (-4)^n
\end{equation}
and the classical expansion
\begin{equation}
(1+x)^{\alpha}=1+\binom{\alpha}{1}x+\binom{\alpha}{2}x^2+\binom{\alpha}{3}x^3+\dots
\end{equation}
we have 
\begin{equation}
\label{eq35}
\sum_{n=0}^{+\infty} u_{2n}=\sum_{n=0}^{+\infty} \binom{2n}{n}p^n q^n=\sum_{n=0}^{+\infty}\binom{-\frac{1}{2}}{n}(-4pq)^n=\frac{1}{\sqrt{1-4pq}}=\frac{1}{|p-q|}
\end{equation}
The last equality is justified by: $(p+q)^2=1\Rightarrow (p-q)^2=1-4pq$.
So we have just proved that, if $p\neq q$, $u:=\sum_{n=0}^{+\infty} u_{2n}<+\infty$, which equals saying $u:=\sum_{n=0}^{+\infty} {\cal P}(E_{2n})<+\infty$.

We now invoke a celebrated theorem, named first Borel-Cantelli Lemma, to draw our conclusion. First Borel-Cantelli Lemma says that if the sum of the probabilities of the events $\sum_{n=0}^{+\infty} {\cal P}(E_{2n})$ is finite, then the probability that infinitely many events $E_{2n}$ occur is $0$, or ${\cal P}( \limsup_{n\to\infty} E_{2n})=0$.

Let us explain the meaning of $\limsup_{n\to\infty} E_{2n}$.
We know that ${\cal P}(E_{2n})={\cal P}(S_{2n}=0)$. If we consider $\cup_{n=1}^{+\infty}E_{2n}$, we are taking all the events $E_{2n}$, for \textit{any} $2n$, from the first trial to the 'last' one, which is the infinite. In other words, we are taking all the possible returns to the origin $E_{2n}$ from time $2n=2$ to infinite. If we consider $\cup_{n=2}^{+\infty}E_{2n}$, we are taking all the events $E_{2n}$ from time $2n=4$ to infinite.
If we consider $\cup_{n=N}^{+\infty}E_{2n}$, we are taking all the events $E_{2n}$ occurred from and after time $2n=2N$. In this way, we do not consider events occurred before time $2n=2N$, that is we are excluding events in the first stage of our experiment. Of course, when $N$ increases, the number of events left decreases and we are restricting to the only events happened after such a time. Now let us imagine that $2N\to +\infty$. There will always be a certain number of events remaining after time 2N, however great that time may be. These events left are in the intersection of all the previous sets, obtained for each finite $N$. The set $\limsup_{n\to\infty} E_{2n}$ is exactly this intersection, the core of events left after each step in which we increase $N$. Formally, $\limsup_{n\to\infty} E_{2n}=\cap_{N=1}^{\infty}\cup_{n=N}^{+\infty}E_{2n}$.

Now the probability related to the core of events left after $2N$ is given by ${\cal P}\left(\cup_{n=N}^{+\infty}E_{2n}\right)$. It is a general law of probability that
${\cal P}\left(\cup_{n=N}^{+\infty}E_{2n}\right)\leq\sum_{n=N}^{+\infty}{\cal P}\left(E_{2n}\right)$.

On the one hand, when $N$ grows to infinity, $\cup_{n=N}^{+\infty}E_{2n}$ becomes our core set, the intersection of all these sets, and so ${\cal P}\left(\cup_{n=N}^{+\infty}E_{2n}\right)$ becomes ${\cal P}( \limsup_{n\to\infty} E_{2n})$. On the other hand, 
when $N$ grows to infinity, $\sum_{n=N}^{+\infty}{\cal P}\left(E_{2n}\right)$ is finite for sure since $\sum_{n=1}^{+\infty}{\cal P}\left(E_{2n}\right)$ is finite, as we proved above. This implies that the limit of the remainder series is $0$, that is $\lim_{N\to +\infty}\sum_{n=N}^{+\infty}{\cal P}\left(E_{2n}\right)=0$. 
Therefore, in the light of the inequality a few lines above, we deduce that also $\lim_{N\to +\infty}{\cal P}\left(\cup_{n=N}^{+\infty}E_{2n}\right)=0$.

This means only one thing: from a certain point forward, the probability of an event $E_{2n}$ occurring is zero, or equivalently, that it can no longer happen that $S_{2n}$ vanishes. The cumulative sums $S_{2n}$ will vanish only \textit{finitely many times}. The random walk will return to zero only in a finite number of times. When this happens, i.e. when, like in this case, we have only a finite number of returns to the origin, we say that the recurrent event $E=$\{A return to the origin takes place at time $2n$\} is {\it transient}.

What does this theorem tell us about gambling?
In gambling terms, it says that after a finite number of initial fluctuations around $0$ the net gain will be positive and remain so if $p>q$ or will be negative and remain so if $p<q$. 

We now merge everything we have said so far, particularly in section \ref{section3}, to demonstrate the probability that the random walk ever re-enters zero. The probability to reach $A$ before $-B$ obtained in Eq. (\ref{eq26}) is
\begin{equation}
{\cal P}(S_{\tau}=A)=\frac{1-\left( \frac{q}{p}\right)^{B}}{1-\left( \frac{q}{p}\right)^{A+B}}
\end{equation}
We can adapt this formula to the case in which $S_0=1$, $B=0$ and $A=N$ (equivalent to $S_0=0$, $B=1$ and $A=N-1$)\footnote{Remind that the lower bound is set to $S_{\tau}=-B$ with $B>0$.}:
\begin{equation}
{\cal P}(S_{\tau}=N|S_0=1)=\frac{1- \frac{q}{p}}{1-\left( \frac{q}{p}\right)^{N}}
\end{equation}

As $N\to +\infty$, the probability that the particle reaches $N$ before $0$ becomes 
\begin{eqnarray*}
\lim_{N\to+\infty }{\cal P}(S_{\tau}=N|S_0=1)=
&&\left\{ 
\begin{array}{c}
0 \qquad \quad \, {\rm if} \ p<q\ {\rm i.e.} \ p<\frac{1}{2}\\ 
1-\frac{q}{p} \  \quad {\rm if} \ p>q \ {\rm i.e.} \ p>\frac{1}{2}\\ 
\end{array}%
\right.
\end{eqnarray*}

Note that when $N\to +\infty$, the upper bound barrier goes to infinity, that is it disappears. The probability that the particle reaches $0$ before $N$ is then
\begin{equation*}
{\cal P}(S_{\tau}\ {\rm reaches\ 0\ before\ N}|S_0=1)= 1-{\cal P}(S_{\tau}\ {\rm reaches\ N\ before\ 0}|S_0=1)
\end{equation*}

and, in the limiting case, it becomes
\begin{eqnarray*}
{\cal P}(S_{\tau}\ {\rm reaches\ 0\ before\ N}|S_0=1)=
&&\left\{ 
\begin{array}{l}
1 \quad \, {\rm if} \ p<q\ {\rm i.e.} \ p<{1}/{2}\\ 
\frac{q}{p}  \quad {\rm if} \ p>q \ {\rm i.e.} \ p>{1}/{2}\\ 
\end{array}%
\right.
\end{eqnarray*}

In a similar way, we can prove that
\begin{eqnarray*}
{\cal P}(S_{\tau}\ {\rm reaches\ 0\ before\ N}|S_0=-1)=
&&\left\{ 
\begin{array}{l}
\frac{p}{q} \quad \, {\rm if} \ p<q\ {\rm i.e.} \ p<{1}/{2}\\ 
1  \quad \ {\rm if} \ p>q \ {\rm i.e.} \ p>{1}/{2}\\ 
\end{array}%
\right.
\end{eqnarray*}

Just a few more steps and we are there.
Let us call ${\cal P}_i={\cal P}(S_{\tau}\ {\rm ever\ returns}\ i\ |S_0=i)$ the probability that the random walk, starting from a point $i$, will sooner or later return to the same point. Then
\begin{equation*}
\begin{split}
{\cal P}_0 &  = {\cal P}(S_{\tau}\ {\rm ever\ returns}\ 0\ |S_0=0)\\
    &  = p\cdot {\cal P}(S_{\tau}\ {\rm ever\ returns}\ 0\ |S_1=1)+q\cdot {\cal P}(S_{\tau}\ {\rm ever\ returns}\ 0\ |S_1=-1)
\end{split}
\end{equation*}
If $p>\frac{1}{2}$ we have
\begin{equation*}
{\cal P}_0 = p\cdot \frac{q}{p} + q\cdot 1= 2q = 1-(p-q)
\end{equation*}
If $p<\frac{1}{2}$ we have
\begin{equation*}
{\cal P}_0 = p\cdot 1 + q\cdot \frac{p}{q}= 2p = 1-(q-p)
\end{equation*}
In general, we have
\begin{equation}
\label{eq38}
{\cal P}_0= 1-|p-q|
\end{equation}

Note that ${\cal P}_0$ represents the probability that the random walk, starting from zero, will sooner or later return to zero. 
We can conclude that

\begin{itemize}
	\item if $p\neq q$ then ${\cal P}_0<1$ and the point $0$ is a transient state,
	\item if $p= q$ then ${\cal P}_0=1$ and the point $0$ is a persistent state,
\end{itemize}

This result is very significant. It implies that, if $p\neq q$, the particle will pass a certain number of times from 0 and then will continue towards constantly positive or negative values according to the value of $p$ and $q$. In other words, the probability that it returns to $0$ more than $r$ times tends to zero. In the case $p=q$, the probability that the particle returns to $0$ more than $r$ times remains equal to $1$ and it will return infinitely many times to $0$.

From all this we conclude that if a game of chance is unfavorable to us, because for us $p$ is less than $q$, then after a certain finite number of fluctuations around zero in the early stages of the game, in which we have now won and now lost, it will eventually be our fate never to break even again, but to continue losing indefinitely!

Let us note that in the previous remarks we also proved the following important relation between ${\cal P}_0$ and $u=\sum_{n=0}^{+\infty} u_{2n}$. If we compare results in Eq. (\ref{eq35}), $u=\frac{1}{|p-q|}$, and in Eq. (\ref{eq38}), ${\cal P}_0= 1-|p-q|$, we find out that
\begin{equation}
u=\frac{1}{1-{\cal P}_0}\qquad {\rm or}\qquad {\cal P}_0=\frac{u-1}{u}
\end{equation} 
Therefore, if $u$ is finite then ${\cal P}_0<1$ and the passage by the origin is transient; if ${\cal P}_0=1$, $u$ diverges and the passage by the origin is persistent. This result will remain valid in dimensions greater than one, as discussed in the Appendix B.

\section{Conclusion}
We conclude by summarizing some of the most important lessons about the unexpected dynamics of gambling that the study of the random walk has revealed so far.

We first observed that, in a fair game, the probability of spending a very small or very large fraction of the time on the positive side, that is in the role of the winner, is much greater than the probability of spending half the time on that side, as our intuition might incorrectly suggest. Moreover, the difference between the probabilities of the extreme cases and that of the $50\%-50\%$ case increases with the number of trials. It does not grow very fast, but it grows! We found that it grows as $\sqrt{N/2}$, where $N$ is the number of trials.

We also learned that there is a probability of $1\%$ that the less fortunate player will be in the lead for $30$ minutes at most, over an entire year of play, while the more fortunate player will be in the lead at least for $364$ days, $23$ hours, and $30$ minutes. There is a non-negligible probability that the unluckiest player will only be in the lead for a very, very small fraction of the year. But there is also a $50\%$ probability that the unlucky player will be in the lead for only $53$ days out of $365$.

The probability of a tie is also very low. The hope of returning to $0$ at least once, which would mean regaining what we lost as less fortunate players, is less than $8$ cases out of $100$ of repeated attempts, and the hope of returning to a tie several times becomes smaller and smaller. And all this in a fair game!

When we moved to the unfair games, things got worse, as we might have expected.
For example, we found that if our strategy is to quit the game when we win A dollars or when we lose our budget of B dollars, then the probability of reaching the win before losing our entire budget becomes very, very small when the game becomes even the slightest bit unfavorable to us. For example, it only takes a $q/p=2$ ratio to be virtually certain of losing \$10 before winning \$10!

Moreover, when again, the probability of loosing is double than the probability of winning and we take $20$ trials, the values of the probability that we break even exactly at the end of the game is $5,43\%$ and the probability that is happens for the first time at the end of the game is only $0,29\%$. As can be seen, the probability of breaking even at the end of the game is very low.

Finally, through a nontrivial argument, we proved that in a unfair game, although it may seem weird, the event related to a return to the origin is {\it transient}, that is it will occur only a finite number of times even if we repeat the trials \textit{ad infinitum}. A return to the origin cannot happen infinitely many times. Only in a fair game, the probability that the random walk, starting from zero, will ever return to zero is equal to $1$ and we are sure that it will, sooner or later, return to the origin. 

\hfill

{\bf Author Declarations}

The author has no conflicts to disclose.
The manuscript has not been published elsewhere and it has not been submitted simultaneously for publication elsewhere.


\appendix

\section{Proofs}

\subsection{Proof of Proposition 1}

Let's take two points $A=(a,y_a)$ and $B=(b,y_b)$ in the positive quadrant with $a<b$ and $y_a, y_b >0$. By {\it reflection} of $A$ on the $x$-axis, we mean the point $A'=(a,-y_a)$. An important lemma links the number of paths from $A$ to $B$ to the number of paths from $A'$ to $B$ and plays a crucial role in the proof of Proposition \ref{theorem1}:

\begin{lemma}
The number of paths from $A$ to $B$ which touch or cross the $x$-axis equals the number of all paths from $A'$ to $B$.
\end{lemma}
\begin{figure}[H]
	\centering
	\includegraphics[width=0.8\textwidth]{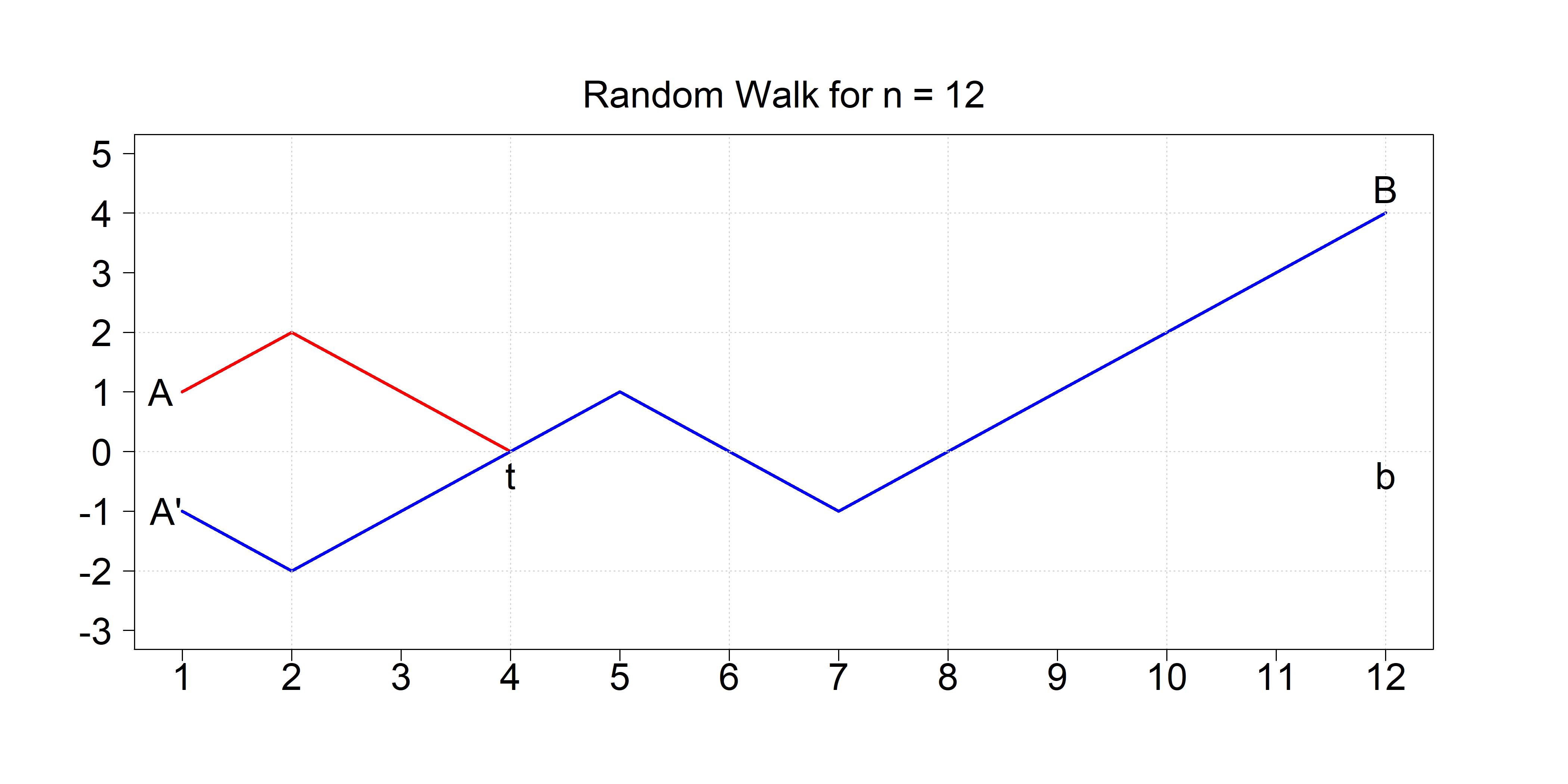}
	\caption{Paths analyzed in the proof of the lemma}
	\label{fig10}
\end{figure}
This preliminary idea is called {\it Reflection Principle}.

Proof of the Lemma: consider a path $y_a=S_a$, $S_{a+1}$, $S_{a+2}$, $\dots$, $S_{b-1}$, $S_{b}=y_b$ from $A$ to $B$ having one or more vertices on the $x$-axis. Let $t$ be the abscissa of the first of such vertices, that is $S_a=y_a>0$, $S_{a+1}>0$, $\dots$, $S_{t-1}>0$, $S_{t}=0$. Then $-S_a=-y_a$, $-S_{a+1}$, $\dots$, $-S_{t-1}$, $S_t=0$, $S_{t+1}$, $S_{t+2}$, \dots, $S_{b-1}$  $S_{b}=y_b$, is a path leading from $A'$ to $B$ and having $T=(t,0)$ as its first point on the $x$-axis. There is so a one-to-one correspondence between paths from $A'$ to $T$ and those from $A$ to $T$, and so between paths from $A'$ to $B$ and those from $A$ to $B$.

Proof of the main Proposition: $S_1$ could be $\pm 1$; but we are looking for positive paths, so it must be $S_1=+1$. Now, we are in the point  $(1,1)$ and we want to find out the number of paths from $(0,0)$ to the point $(x,y)$ which neither touch nor cross the $x$-axis. By previous lemma, the number of such paths is equal to
\begin{equation*}
N_{x-1,y-1}-N_{x-1,y+1}
\end{equation*}
Here it is why. Let's have a look at the Fig. \ref{fig11}:
\begin{figure}[H]
	\centering
	\includegraphics[width=0.8\textwidth]{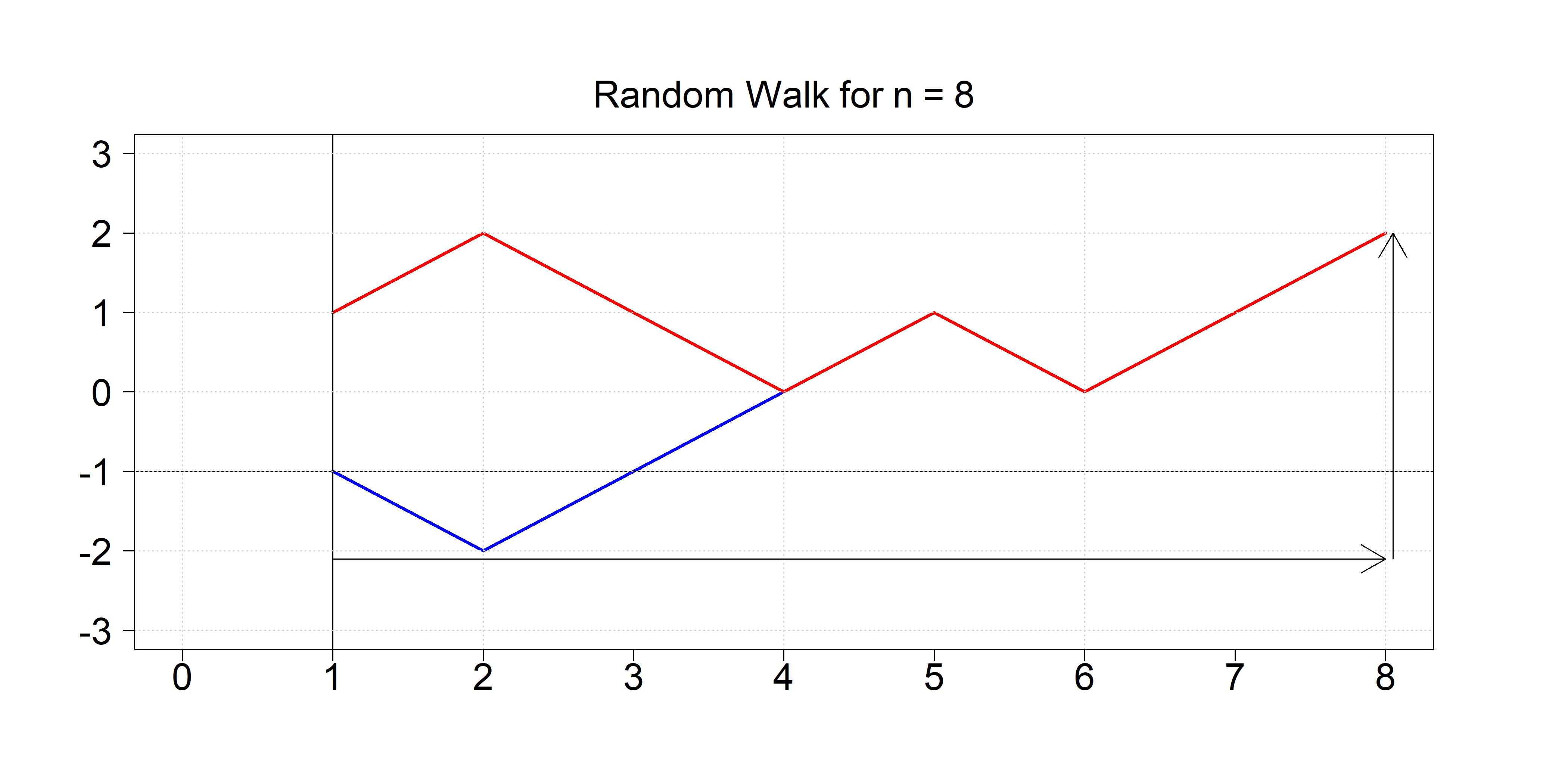}
	\caption{Paths analyzed in the proof of the Proposition}
	\label{fig11}
\end{figure}

We have
\begin{itemize}
\item $N_{x,y}$ is the number of all possible paths from $(0,0)$ to $(x,y)$;
\item $N_{x-1,y-1}$ is the number of all possible paths from $(0,0)$ to $(x-1,y-1)$, which is the same as the number of paths from $(1,1)$ to $(x,y)$;
\item $N_{x-1,y+1}$ is the number of all possible paths from $(0,0)$ to $(x-1,y+1)$, which is the same as the number of paths from $(1,-1)$ to $(x,y)$.
\end{itemize}
The set of all the paths from from $(1,1)$ to $(x,y)$ that do touch the $x$-axis anywhere is in a one-to-one correspondence with the set of all paths from $(1,-1)$ to $(x,y)$, since they have to cross the $x$-axis in at least one point. This is established by reflection with respect to the $x$-axis of the initial segment of the path at the step where it first touches the $x$-axis. This number is given by $N_{x-1,y+1}$ and it can be seen if we change the coordinate axes moving the origin to $(-1,1)$. So, now we have to subtract from the number of all paths from $(1,1)$ to $(x,y)$ the number of paths that start in $(1,1)$ and end in $(x,y)$ and that touches the $x$-axis in at least one point. Coming back to our language in terms of $p$ and $q$ and being $x=p+q$ and $y=p-q$, we have
\begin{equation*}
N_{x,y}=\binom{p+q}{p} \quad N_{x-1,y-1}=\binom{p+q-1}{p-1} \quad N_{x-1,y+1}=\binom{p+q-1}{q-1}
\end{equation*}
from which
\begin{equation*}
\begin{split}
&N_{x-1,y-1}-N_{x-1,y+1}\\
=&\binom{p+q-1}{p-1} -\binom{p+q-1}{q-1}\\
=&\frac{p-q}{p+q}\binom{p+q}{p}=\frac{y}{x}N_{x,y}
\end{split}
\end{equation*}
and the proposition is proved. 

\subsection{Proof of Proposition 2}

Each path such that $S_{1}> 0$, $S_{2}> 0$, $\dots$, $S_{2n-1}> 0$ and $S_{2n}= 0$ must pass through the point $(2n-1,1)$ and the number of paths such that $S_{1}> 0$, $S_{2}> 0$, $\dots$, $S_{2n-2}> 0$ equals
\begin{equation*}
\frac{y}{x}N_{x,y}=\frac{1}{2n-1}\binom{2n-1}{n-1}=\frac{1}{n}\binom{2n-2}{n-1}
\end{equation*}
since $x=p+q=2n-1$, $y=p-q=1$ and so $p=n$ and $q=n-1$. This is equal to $L_{2n-2}$ and it proves the first statement. Now let us consider a path joining $(1,1)$ and $(2n-1,1)$: if $S_{1}> 0$, $S_{2}> 0$, $\dots$, $S_{2n-1}> 0$, then all its vertices lie on or above the line $y=1$. Now, translating the origin to $(1,1)$, these paths connect the point $(0,0)$ to $(2n-2,0)$ and all their vertices are on or above the $x$-axis. We have established a one-to-one correspondence between the paths satisfying $S_{1}\geq 0$, $S_{2}\geq 0$, $\dots$, $S_{2n-1}\geq 0$ and those such that $S_{1}> 0$, $S_{2}> 0$, $\dots$, $S_{2n-1}> 0$ but with $n$ replaced by $n+1$ because we moved the origin to $(1,1)$. So this number is $L_{2(n+1)-2}=L_{2n}$.

\subsection{Proof of Proposition 3}
Let us consider the space of paths of fixed length $2n$. The total number of paths from $(0,0)$ to $(2n,y),\ \forall y$, is $2^{2n}$. We proved there exist
\begin{equation*}
N_{2n,0}=\binom{2n}{n}
\end{equation*}
paths from $(0,0)$ to $(2n,0)$, so the relation a) is proved. We have seen before that there exist $L_{2n-2}$ paths joining $(0,0)$ to $(2n,0)$ such that $S_{1}> 0$, $S_{2}> 0$, $\dots$, $S_{2n-1}> 0$. Therefore there are twice as many paths such that $S_{1}\neq 0,\dots , S_{2n-1} \neq 0$ and the corresponding probability is $2L_{2n-2}2^{-2n}=2^{-2n+1}L_{2n-2}$ and this gives relation $d)$. The number of paths such that $S_{1}\geq 0,\dots , S_{2n-2} \geq 0, S_{2n-1} < 0$ is equal to the number of paths such that $S_{1}\geq 0,\dots , S_{2n-3} \geq 0, S_{2n-2} = 0$ which is given by $L_{2(n-1)}=L_{2n-2}$ and again we obtain probability dividing by the total number of paths that, in this case, is $2\cdot 2^{2(n-1)}=2^{2n-1}$ because after $2n-2$ we could go up or down. So again we get $2^{-2n+1}L_{2n-2}$ and $e)$ is proved. The probability that no zero occurs up to and including time $2n$ is complementary to the probability that there is a first return to the origin at any of the even times less than or equal $2n$. So we have $1-(u_0-u_2)-(u_2-u_4)-\cdots -(u_{2n-2}-u_{2n})= u_{2n}$ and this proves $b)$. Finally, the probability that we have paths such that $S_{1}\geq 0,\dots , S_{2n} \geq 0$ is complementary to the probability of paths with a first passage to $-1$ at any time less than or equal to $2n$, and, as before, we get $u_{2n}$. Finally, to prove the last equality in $d)$ and $e)$ we have to take into account the combinatorial relation
$
u_{2n-2}=\frac{2n}{2n-1}u_{2n}.
$

\subsection{Proof of Proposition 4}
Let us observe at first that, if $2k=0$ or $2k=2n$, we have $p_{0,2n}=p_{2n,2n}=u_{2n}$ according to point c) of proposition \ref{theorem3}. Let $1\leq k\leq n-1$. A particle that stays on the positive side for $2k>0$ time units and on the negative side for $2n-2k>0$ time units must pass through zero. Let $2t$ be the time of its first return to zero. We can have two cases depending on whether the particle was on the positive or negative side before $2t$.

1) Up to time $2t$ the particle stays on the positive side and during the interval $(2t,2n)$ it spends exactly $2k-2t$ time units on the positive side:
\begin{figure}[H]
	\centering
	\includegraphics[width=0.8\textwidth]{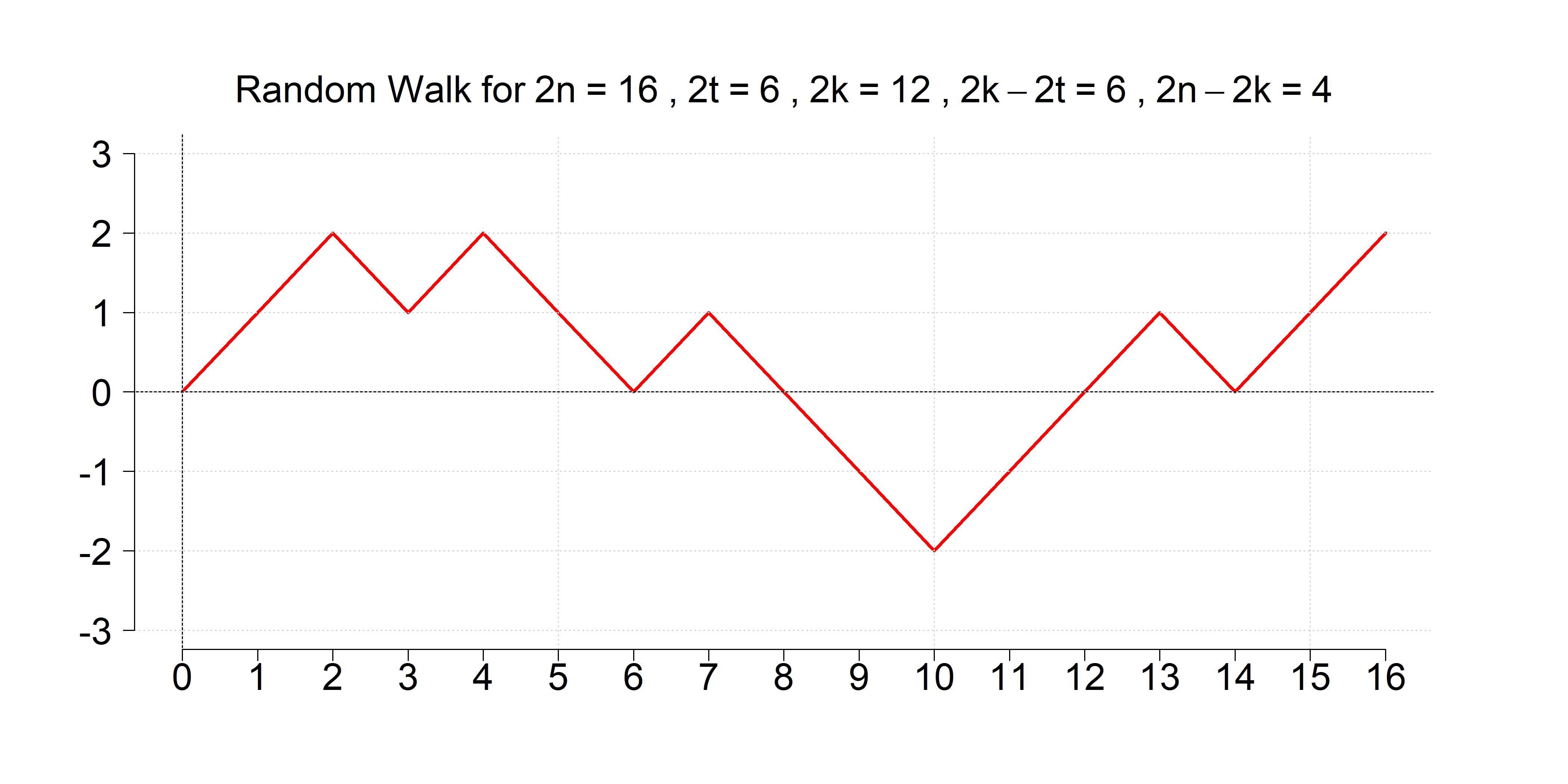}
	\caption{Positive path from $0$ to $2t$.}
	\label{fig12}
\end{figure}
Then there exist $2^{2t}\cdot (u_{2t-2}-u_{2t})$ paths of length $2t$ which return to the origin for the first time at $2t$. It is the total number of paths $2^{2t}$ multiplied by the probability $u_{2t-2}-u_{2t}$ given by point d) of proposition \ref{theorem3}. Half of these paths keeps to the positive side: $2^{2t-1}\cdot (u_{2t-2}-u_{2t})$. Furthermore, by definition, there are $2^{2n-2t}\cdot p_{2k-2t,2n-2t}$ paths of length $2n-2t$ starting at $(2t,0)$ and having exactly $2k-2t$ sides above the $x$-axis. Thus the total number of paths of length $2n$ of the first type is
\begin{equation*}
2^{2t-1}(u_{2t-2}-u_{2t})\cdot 2^{2n-2t}\cdot p_{2k-2t,2n-2t}=2^{2n-1} (u_{2t-2}-u_{2t})p_{2k-2t,2n-2t}
\end{equation*}
2) From $0$ to $2t$ the particle stays on the negative side, and between $2t$ and $2n$ it spends exactly $2k$ time units on the positive side:
\begin{figure}[H]
	\centering
	\includegraphics[width=0.8\textwidth]{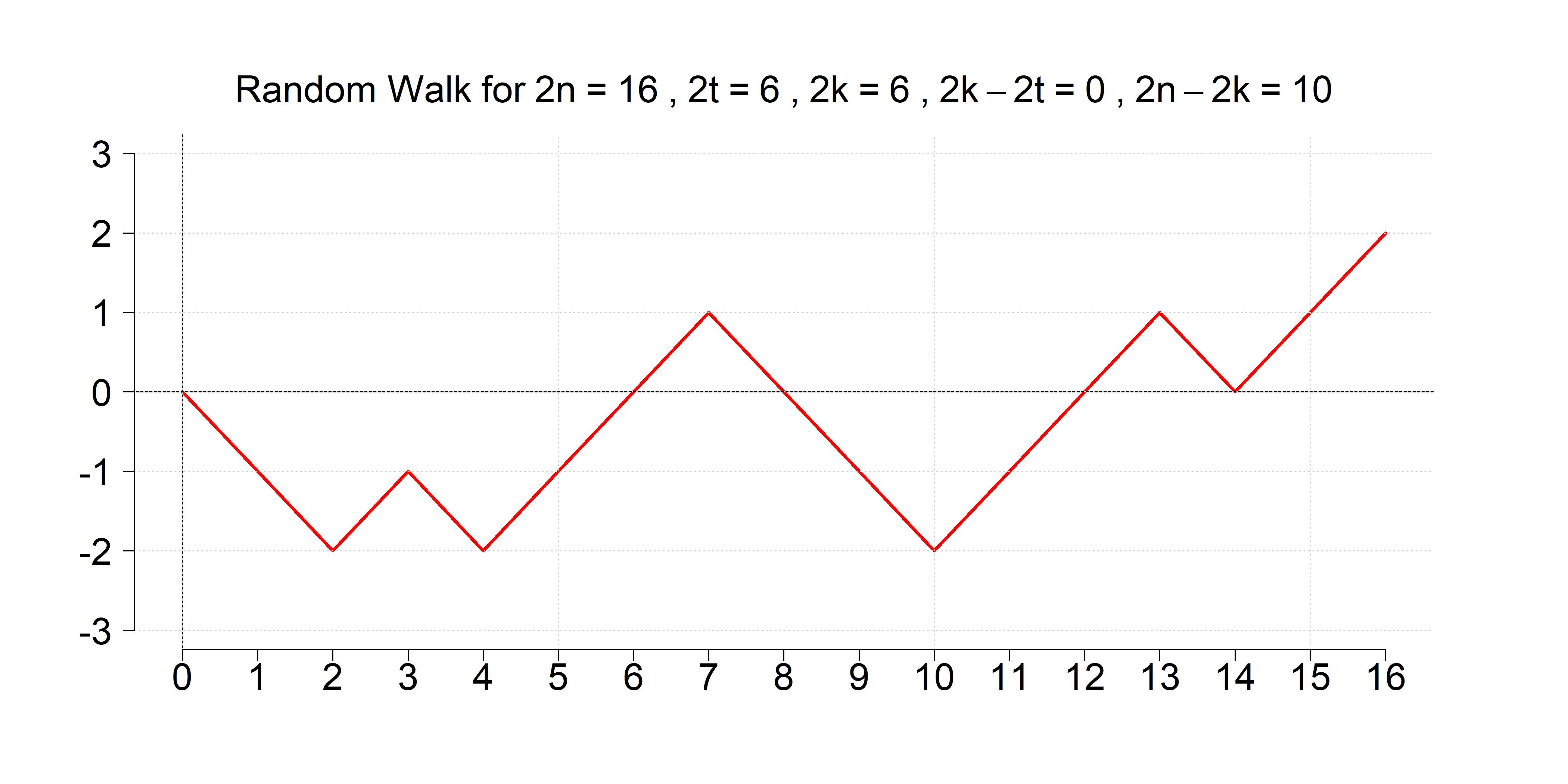}
	\caption{Negative path from $0$ to $2t$.}
	\label{fig13}
\end{figure}
Of course, now $2k+2t\leq 2n$, i.e. $2k\leq 2n-2t$. All the $2k$ positive intervals are at the right  of $2t$, after $2t$. So we can use the same argument as before but with $2k$ instead of $2k-2t$ and the number of paths is 
\begin{equation*}
2^{2n-1} (u_{2t-2}-u_{2t})p_{2k,2n-2t}
\end{equation*}
It follows that, for $1\leq k\leq n-1$:
\begin{equation*}
\small
\begin{split}
& p_{2k,2n}=\\
& \ \frac{1}{2^{2n}}\left[ \sum_{t=1}^{k} 2^{2n-1} (u_{2t-2}-u_{2t})p_{2k-2t,2n-2t}+\sum_{t=1}^{n-k} 2^{2n-1} (u_{2t-2}-u_{2t})p_{2k,2n-2t} \right]=  \\
& \ \frac{1}{2} \sum_{t=1}^{k} (u_{2t-2}-u_{2t})p_{2k-2t,2n-2t}+ \frac{1}{2} \sum_{t=1}^{n-k} (u_{2t-2}-u_{2t})p_{2k,2n-2t}=   \\
& \ \frac{1}{2} \left[ (u_0-u_2)p_{2k-2,2n-2}+(u_2-u_4)p_{2k-4,2n-4}+\dots + (u_{2k-2}-u_{2k})p_{0,2n-2k} \right]+\\
&\ \frac{1}{2} \left[ (u_0-u_2)p_{2k,2n-2}+(u_2-u_4)p_{2k,2n-4}+\dots + (u_{2n-2k-2}-u_{2n-2k})p_{2k,2k} \right] \\
\end{split}
\end{equation*}

\noindent Now let us suppose by induction that $p_{2k,2\nu}=u_{2k}u_{2\nu-2k}$ for $\nu=1,2,\dots, n-1$. Then we have:
\begin{equation*}
\small
\begin{split}
& p_{2k,2n}=\\
& \frac{1}{2} [ (u_0-u_2)u_{2k-2}u_{2n-2k}+(u_2-u_4)u_{2k-4}u_{2n-2k}+ \dots + (u_{2k-2}-u_{2k})u_{0}u_{2n-2k} + \\
&(u_0-u_2)u_{2k}u_{2n-2k-2}+(u_2-u_4)u_{2k}u_{2n-2k-4}+ \dots + (u_{2n-2k-2}-u_{2n-2k})u_{2k}u_{0}]= \\
& \frac{1}{2} u_{2n-2k} \left( u_0 u_{2k-2}-u_2 u_{2k-2}+ u_2 u_{2k-4}-u_4 u_{2k-4}+ \dots + u_{2k-2} u_{0}-u_{2k}u_{0} \right) + \\
&\frac{1}{2} u_{2k} \left( u_0 u_{2n-2k-2}-u_2 u_{2n-2k-2}+ u_2 u_{2n-2k-4}-u_4 u_{2n-2k-4}+ \dots + u_{2n-2k-2} u_{0}-u_{2n-2k}u_{0} \right) \\
=& \frac{1}{2} u_{2n-2k} u_{2k}+ \frac{1}{2} u_{2k} u_{2n-2k} = u_{2k} u_{2n-2k}.
\end{split}
\end{equation*}

\subsection{Proof of Proposition 5}
Let us observe that, by means of the Stirling's approximation formula
\begin{equation*}
u_{2n}=\frac{1}{2^{2n}}\binom{2n}{n}=\frac{1}{2^{2n}}\frac{(2n)!}{(n!)^2}\sim \frac{\sqrt{2\pi 2n}\left(\frac{2n}{e}\right)^{2n}}{2\pi n \left(\frac{n}{e}\right)^{2n}2^{2n}}=\frac{1}{\sqrt{\pi n}}
\end{equation*}
So that, if we define the ratio $\frac{2k}{2n}$ as $\rho=\frac{k}{n}$
\begin{equation*}
p_{2k,2n}=u_{2k}\cdot u_{2n-2k} \sim \frac{1}{\sqrt{\pi k}}\cdot \frac{1}{\sqrt{\pi (n-k)}}=\frac{1}{\pi n \sqrt{\frac{k}{n}\left( 1-\frac{k}{n}\right)}}=\frac{1}{\pi n \sqrt{\rho( 1-\rho)}}
\end{equation*}
and
\begin{equation*}
{\cal P}(2k\leq 2\alpha) =\sum_{\rho=0}^{\frac{\alpha}{n}} \frac{1}{\pi n}\frac{1}{\sqrt{\rho( 1-\rho)}}
\end{equation*}
For $n, k, \alpha \to +\infty$, we recognize the Riemann sum approximating the integral
\begin{equation*}
\int_{0}^{\alpha}\frac{1}{\pi} \frac{dk}{n}\frac{1}{\sqrt{\frac{k}{n}\left( 1-\frac{k}{n}\right)}}= \frac{1}{\pi} \int_{0}^{\alpha /n}\frac{d\rho}{\sqrt{\rho( 1-\rho)}}
\end{equation*}
By setting $\rho=\cos^2\theta$ and $\theta = {\rm arccos}\, \sqrt{\rho}$, we get
\begin{equation*}
{\cal P}(2k\leq 2\alpha) = \frac{1}{\pi} \left[ -2{\rm arccos}\, \sqrt{\rho} \right]_{0}^{\alpha /n} =\frac{2}{\pi}\left[ \frac{\pi}{2} -{\rm arccos}\, \sqrt{\frac{\alpha}{n}} \right]
\end{equation*}
and eventually we have
\begin{equation*}
{\cal P}(2k\leq 2\alpha)= \frac{2}{\pi}{\rm arcsin}\, \sqrt{\frac{\alpha}{n}}
\end{equation*}

\subsection{Proof of Proposition 6}
In probability theory, a sequence of random variables for which, at any given time, the conditional expectation of the subsequent value is equal to the current one, regardless of any previous value, is called a martingale. This proof will use some properties of martingales, and we refer the reader to the text \cite{Williams} for some of the technical details, particularly in reference to the stopped process and the stopping time theorem.

The process $S_n$ is a martingale and, by the stopping time theorem, also $S_{n\wedge\tau}$, the stopped process, is a martingale. It follows that ${\mathbb E}[S_{n\wedge\tau}]={\mathbb E}[S_{0\wedge\tau}]={\mathbb E}[S_{0}]=0$ if $S_0=0$. So ${\mathbb E}[S_{n\wedge\tau}]=0,\ \forall n \geq 0$. We also ask that ${\cal P}(\tau<+\infty)=1$, i.e. that $\tau$ is finite with probability $1$. Under this conditions: $\lim_{n\to +\infty} S_{n\wedge\tau}=S_{\tau}$, almost surely. Since $|S_{n\wedge\tau}|\leq \max (A,B)$ we can apply the dominated convergence theorem to conclude that
\begin{equation*}
\lim_{n\to +\infty} {\mathbb E}[S_{n\wedge\tau}]={\mathbb E}[S_{\tau}]
\end{equation*}
So we have also ${\mathbb E}[S_{\tau}]=0$, but ${\mathbb E}[S_{\tau}]=A\cdot {\cal P}(S_{\tau}=A)-B\cdot {\cal P}(S_{\tau}=-B)$. Therefore, it follows that
\begin{equation*}
A\cdot {\cal P}(S_{\tau}=A)-B\cdot \left( 1-{\cal P}(S_{\tau}=A)\right)=0 \ \Rightarrow (A+B){\cal P}(S_{\tau}=A)=B,
\end{equation*}
from which we have the thesis. In order to prove the second part, we observe that ${\rm Var}[X_n]={\mathbb E}[X_n^2]-{\mathbb E}^2[X_n]={\mathbb E}[X_n^2]=1$ so $M_n :=S_n^2-n$ is also a martingale. Reasoning as before, ${\mathbb E}[M_\tau]=0$ but
\begin{equation*}
{\mathbb E}[M_\tau]={\mathbb E}[S_{\tau}^2]-{\mathbb E}[\tau]=A^2 \cdot {\cal P}(S_{\tau}=A)+B^2 \cdot {\cal P}(S_{\tau}=-B)-{\mathbb E}[\tau]
\end{equation*}
so that
\begin{equation*}
{\mathbb E}[\tau]=A^2\cdot \frac{B}{A+B}+B^2\cdot \frac{A}{A+B}=AB.
\end{equation*}

\subsection{Proof of Proposition 7}
Let us define a new process $M_n:=\big(\frac{q}{p}\big)^{S_n}$, with $M_0=1$. The process $M_n$ is a martingale, since
\begin{equation*}
\begin{split}
&{\mathbb E}[M_{n+1}|{\cal F}_{n}]={\mathbb E}\left[\left(\frac{q}{p}\right)^{S_{n+1}}|{\cal F}_{n}\right]={\mathbb E}\left[\left(\frac{q}{p}\right)^{S_{n}+X_{n+1}}|{\cal F}_{n}\right]=\\
&{\mathbb E}\left[\left(\frac{q}{p}\right)^{S_{n}}\cdot \left(\frac{q}{p}\right)^{X_{n+1}}|{\cal F}_{n}\right]= \left(\frac{q}{p}\right)^{S_{n}}\cdot {\mathbb E}\left[\left(\frac{q}{p}\right)^{X_{n+1}}|{\cal F}_{n}\right]=\\
&\left(\frac{q}{p}\right)^{S_{n}}\cdot\left[ \left(\frac{q}{p}\right)^{1}\cdot {\cal P}(X_n=1) + \left(\frac{q}{p}\right)^{-1}\cdot {\cal P}(X_n=-1)\right]=\\
&\left(\frac{q}{p}\right)^{S_{n}}\cdot\left[\frac{q}{p}p+\frac{p}{q}q\right]= \left(\frac{q}{p}\right)^{S_{n}}\cdot[q+p]= \left(\frac{q}{p}\right)^{S_{n}}=M_n
\end{split}
\end{equation*}
This fact implies that ${\mathbb E}[M_{\tau}]=\lim_{n\to +\infty}{\mathbb E}[M_{n\wedge\tau}]={\mathbb E}[M_0]=1$, but
\begin{equation*}
{\mathbb E}[M_{\tau}]=\left( \frac{q}{p}\right)^{A}\cdot {\cal P}(S_{\tau}=A)+\left( \frac{q}{p}\right)^{-B}\cdot {\cal P}(S_{\tau}=-B)=1
\end{equation*}
so that
\begin{equation*}
\left( \frac{q}{p}\right)^{A}\cdot {\cal P}(S_{\tau}=A)+\left( \frac{q}{p}\right)^{-B}\cdot \left[1- {\cal P}(S_{\tau}=A)\right] =1
\end{equation*}
which equals
\begin{equation*}
\left[ \left( \frac{q}{p}\right)^{A} -\left( \frac{q}{p}\right)^{-B} \right]\cdot {\cal P}(S_{\tau}=A)=1-\left( \frac{q}{p}\right)^{-B}.
\end{equation*}
Finally, we get
\begin{equation*}
{\cal P}(S_{\tau}=A)=\frac{1-\left( \frac{q}{p}\right)^{-B}}{\left( \frac{q}{p}\right)^{A}-\left( \frac{q}{p}\right)^{-B}}=\frac{1-\left( \frac{q}{p}\right)^{B}}{1-\left( \frac{q}{p}\right)^{A+B}}
\end{equation*}
We have also
\begin{equation*}
{\cal P}(S_{\tau}=B)=1-{\cal P}(S_{\tau}=A)=1-\frac{1-\left( \frac{q}{p}\right)^{B}}{1-\left( \frac{q}{p}\right)^{A+B}}=\frac{1-\left( \frac{p}{q}\right)^{A}}{1-\left( \frac{p}{q}\right)^{A+B}}.
\end{equation*}
The last expression is equivalent to ${\cal P}(S_{\tau}=A)$ under the switch $q\leftrightarrow p$ and $A \leftrightarrow B$.
To prove the expected duration, let us consider a new process:
\begin{equation*}
M_n=S_n-{\mathbb E}[X_n]\cdot n
\end{equation*}
In the biased random walk we have ${\mathbb E}[X_n]=1\cdot p + (-1)\cdot q= p-q=2p-1$. So we define: 
\begin{equation*}
M_n=S_n-(p-q)\cdot n
\end{equation*}
$M_n$ is a martingale since
\begin{equation*}
\begin{split}
{\mathbb E}[M_{n+1}|{{\cal F}_n}]=& {\mathbb E}[S_{n+1}-(p-q)(n+1)|{{\cal F}_n}]\\
=&{\mathbb E}[S_n+X_{n+1}-n(p-q)-(p-q)|{{\cal F}_n}]\\
=&S_n-n(p-q)+E[X_{n+1}]-(p-q)\\
=&S_n-n(p-q)+(p-q)-(p-q)\\
=&S_n-n(p-q)=M_n
\end{split}
\end{equation*}
and we can apply stopping theorem for martingales.
In fact we know that ${\mathbb E}[S_{n\wedge\tau}]={\mathbb E}[S_0]=0$ and $\lim_{n\to +\infty}{\mathbb E}[S_{n\wedge\tau}]={\mathbb E}[S_{\tau}]$ because of the theorem of dominated convergence. So we have ${\mathbb E}[M_{\tau}]=0$ and
\begin{equation*}
{\mathbb E}[M_{\tau}]={\mathbb E}[S_{\tau}]-(p-q){\mathbb E}[{\tau}] \Rightarrow (p-q){\mathbb E}[{\tau}]={\mathbb E}[S_{\tau}] \Rightarrow {\mathbb E}[{\tau}]=\frac{1}{p-q}{\mathbb E}[S_{\tau}]
\end{equation*}
Therefore:
\begin{equation*}
\begin{split}
{\mathbb E}[{\tau}]=&\frac{1}{p-q}\left[ A\cdot \frac{1-\big(\frac{q}{p}\big)^{B}}{1-\big(\frac{q}{p}\big)^{A+B}}- B\cdot \left( 1-\frac{1-\big(\frac{q}{p}\big)^{B}}{1-\big(\frac{q}{p}\big)^{A+B}} \right) \right]\\
=&\frac{1}{p-q}\left[ (A+B)\frac{1-\big(\frac{q}{p}\big)^{B}}{1-\big(\frac{q}{p}\big)^{A+B}}-B \right]]\\
=&\frac{A+B}{p-q}\cdot \frac{1-\big(\frac{q}{p}\big)^{B}}{1-\big(\frac{q}{p}\big)^{A+B}}-\frac{B}{p-q}
\end{split}
\end{equation*}
The ${\mathbb E}[\tau]$ may be given a more symmetric expression
\begin{equation*}
{\mathbb E}[\tau]=\frac{1}{p-q}\left[ A \frac{1-\big(\frac{q}{p}\big)^{B}}{1-\big(\frac{q}{p}\big)^{A+B}}- B\frac{1-\big(\frac{p}{q}\big)^{A}}{1-\big(\frac{p}{q}\big)^{A+B}}\right]
\end{equation*}
or
\begin{equation*}
{\mathbb E}[\tau]=\frac{1}{p-q}\frac{A\left( 1-\big(\frac{p}{q}\big)^{B} \right)+B \left( 1-\big(\frac{q}{p}\big)^{A} \right) }{\big(\frac{q}{p}\big)^{A}-\big(\frac{p}{q}\big)^{B}}
\end{equation*}

\subsection{Proof of the Remark after Proposition \ref{theorem7}}
If $p=q$ then $p\rightarrow\frac{1}{2}$ and $q\rightarrow\frac{1}{2}$, and we can consider an expansion of our functions around $\frac{1}{2}$. If we set
\begin{equation*}
f(p)=\left(\frac{q}{p}\right)^A=\left(\frac{1}{p}-1\right)^A
\end{equation*}
we have that 
\begin{equation*}
f'\left( p \right)=-\frac{A}{p^2}\left(\frac{1}{p}-1\right)^{A-1} \Rightarrow \ f'\left( \frac{1}{2} \right)=-4A
\end{equation*}
and
\begin{equation*}
f''\left( p \right)=\frac{2A}{p^3}\left(\frac{1}{p}-1\right)^{A-1}+\frac{A(A-1)}{p^4}\left(\frac{1}{p}-1\right)^{A-2}  \Rightarrow \ f'\left( \frac{1}{2} \right)=16A^2
\end{equation*}
so that
\begin{equation*}
f(p)\sim 1-4A\left(p-\frac{1}{2}\right)+8A^2\left(p-\frac{1}{2}\right)^2
\end{equation*}
which is equivalent to
\begin{equation*}
1-\left( \frac{q}{p}\right) ^A \sim 2A(p-q)-2A^2(p-q)^2
\end{equation*}
Thus we have
\begin{equation*}
{\cal P}(S_\tau=A)=\frac{1-\big(\frac{q}{p}\big)^{B}}{1-\big(\frac{q}{p}\big)^{A+B}} \sim \frac{2B(2p-1)}{2(A+B)(2p-1)}=\frac{B}{A+B}
\end{equation*}
and
\begin{equation*}
{\mathbb E}[\tau]\sim \frac{1}{p-q}\left[ A \frac{2B(p-q)-2B^2(p-q)^2}{2(A+B)(p-q)}- B \frac{2A(q-p)-2A^2(q-p)^2}{2(A+B)(q-p)} \right]
\end{equation*}

\begin{equation*}
=\frac{AB}{p-q}\left[ \frac{B(p-q)+A(p-q)}{A+B}\right]=AB
\end{equation*}

\section{Higher dimensions}

Let us consider now a two-dimensional symmetric random walk that starts at the origin and is performed on the lattice ${\mathbb Z}^2$. We can imagine a walker that now moves on integer points in two dimensions: each step is a distance $1$ jump in one of the four directions - let's say {\it up}, {\it down}, {\it right}, {\it left}. So there are four directions for each step and the choice of a direction is random with probability $\frac{1}{4}$ for each one. We will consider only finite walks. A given walk of length $N$ is then performed with probability $\frac{1}{4^N}$. We define {\it loop} a walk that begins and ends at the origin. A walk of length zero is a {\it trivial} loop and a loop is said to be {\it simple} if it is not a concatenation of two nontrivial loops.


Obviously, a loop has an even length. Let $N^{(2)}_{2n}$ be the number of loops of $2n$ steps. The probability of a return to the origin in $2n$ steps ($2n$ includes all the steps in the $x$- and $y$-directions) is then
\begin{equation*}
{\cal P}\left( S^{(2)}_{2n}=(0,0)\right) :=u^{(2)}_{2n}=\frac{1}{4^{2n}}N^{(2)}_{2n}
\end{equation*}
A return to the origin is possible only if the number of steps  in the positive $x$- and $y$-directions equal those in the negative $x$- and $y$-directions, respectively. So if we divide the $2n$ steps into the four classes {\it up}, {\it down}, {\it right}, {\it left} and if there are $k$ steps in the first one, then we must have $k$ steps in the second one and $n-k$ steps in the third and forth one. The number of all the combinations of these steps is given by a multinomial coefficient and if we sum over all $k$, $0\leq k\leq n$, we get the number of all $2n$-loops. In this way, we have
\begin{equation*}
u^{(2)}_{2n}=\frac{1}{4^{2n}}\sum_{k=0}^{n}\frac{(2n)!}{k!k!(n-k)!(n-k)!}=\frac{1}{4^{2n}}\binom{2n}{n}\sum_{k=0}^{n}{\binom{n}{k}}^2=\frac{1}{4^{2n}}{\binom{2n}{n}}^2
\end{equation*}
since it holds ${\binom{n}{0}}^2+{\binom{n}{1}}^2+\dots +{\binom{n}{n}}^2=\binom{2n}{n}$.

Let us observe that we could have obtained this number also in the following way.
Let us consider two strings $a$ and $b$ of $+1$ and $-1$ of length $2n$, such that there is an equal number of $+1$ and $-1$ in each one. We want to construct a loop of length $2n$. We can establish the following rules: the pair $(a_{i},b_{i})$ corresponds to the $i$-th step

\begin{itemize}
  \item in the direction $(+1,0)$ if $(a_{i},b_{i})= (+1,+1)$,
	\item in the direction $(-1,0)$ if $(a_{i},b_{i})= (-1,-1)$. 
	\item in the direction $(0,+1)$ if $(a_{i},b_{i})= (-1,+1)$,
	\item in the direction $(0,-1)$ if $(a_{i},b_{i})= (+1,-1)$,
\end{itemize}

This walk returns to the origin. In fact, let the number of pairs $(+1,+1)$ be $k$, the number of pairs $(-1,-1)$ be $l$, the number of pairs $(-1,+1)$ be $m$ and the number of pairs $(+1,-1)$ be $n$, then $k+m=l+n$ and $k+n=l+m$. Hence, $2k=2l$, i.e. $k=l$ and $m=n$. So, the number of {\it up} steps is equal to the number of {\it down} steps and the number of {\it right} steps is equal to the number of {\it left} steps.

\begin{figure}[H]
	\centering
	\includegraphics[width=0.6\textwidth]{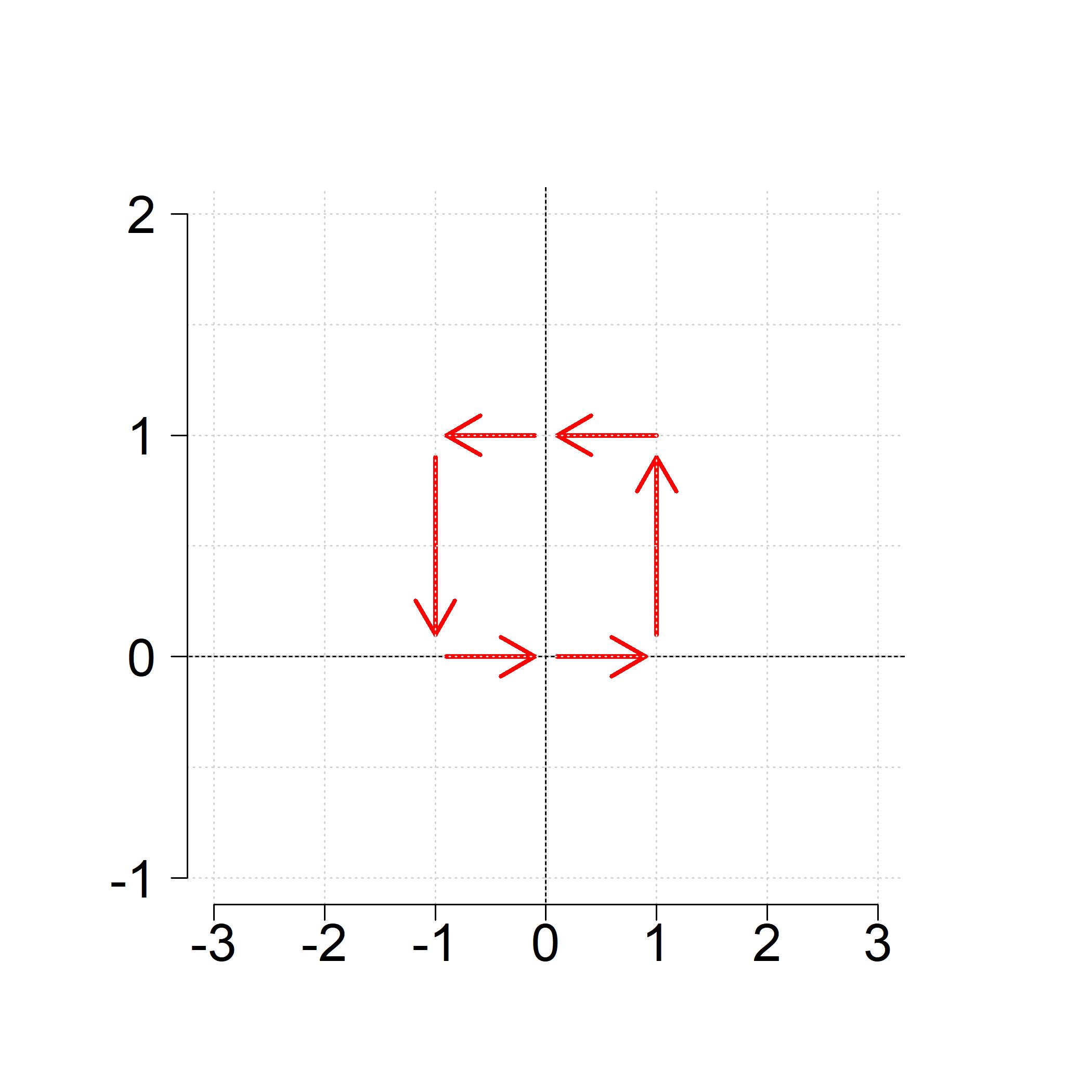}
	\caption{Loop in two dimensions corresponding to the sequence $\binom{+---++}{++---+}$.}
	\label{fig14}
\end{figure}
In this way we have built a bijection between the loops in ${\mathbb Z}^2$ and all the couples of plus and minus that we can build starting from two strings of $2n$ terms in which we have $n$ plus and $n$ minus, and this number is of course the product of the two binomials $\binom{2n}{n}$.
Now, again for $n\to+\infty$ by Stirling's approximation ${\binom{2n}{n}}^2 \sim \frac{4^{2n}}{\pi n}$, we have
\begin{equation*}
u^{(2)}_{2n} \sim \frac{1}{\pi n}
\end{equation*}
Again
\begin{equation*}
u^{(2)}:=\sum_{n=0}^{+\infty} u^{(2)}_{2n} \sim \sum_{n=0}^{+\infty} \frac{1}{\pi n}
\end{equation*}
diverges, so that
\begin{equation*}
{\cal P}^{(2)}_{0}={\cal P}\left( S_{\tau}\ {\rm ever\ returns\ to\ } (0,0)|S_{0}=(0,0)\right)=\frac{u^{(2)}-1}{u^{(2)}}=1
\end{equation*}
This means that there is a probability $1$ that the particle will sooner or later - and therefore infinitely often - return to its initial position, i.e. the particle will pass through every possible point infinitely often.

In the case of three dimensions, our particle can move to six different neighborhoods with equal probability. Thus, if we define
\begin{equation*}
{\cal P}\left( S^{(3)}_{2n}=(0,0,0)\right) :=u^{(3)}_{2n}
\end{equation*}
we have
\begin{equation*}
u^{(3)}_{2n}=\frac{1}{6^{2n}}\sum_{j=0}^{n-k}\sum_{k=0}^{n}\frac{(2n)!}{j!j!k!k!(n-j-k)!(n-j-k)!}
\end{equation*}
the sums being extended over all $j$ and $k$ with $j+k\leq n$. It is easy to see that it is equivalent to 
\begin{equation*}
u^{(3)}_{2n}=\frac{1}{2^{2n}}\binom{2n}{n}\sum_{j,k}\left\{ \frac{1}{3^n}\frac{n!}{j!k!(n-j-k)!}\right\}^2
\end{equation*}
Let us observe that the quantity $p_{j,k}=\frac{1}{3^n}\frac{n!}{j!k!(n-j-k)!}$ represents a trinomial distribution, so that $\sum_{j,k}p_{j,k}=1$. Now we know that if there are $N$ positive numbers $a_1<a_2<\dots <a_N$ between $0$ and $1$ such that $\sum_{i=1}^{N}a_i =1$ then $\sum_{i=1}^{N}a^2_i <\sum_{i=1}^{N}a_i a_N =a_N \sum_{i=1}^{N}a_i=a_N$, i.e. the sum of the squared terms is less than the greatest term. So we can say that the sum of the squares of the $p_{j,k}$'s is less than the biggest one among them and the biggest one is obtained when $j=k=\frac{n}{3}$ since it is a multinomial distribution. From this fact, it follows
\begin{equation*}
u^{(3)}_{2n}<\frac{1}{2^{2n}}\binom{2n}{n} \frac{1}{3^n}\frac{n!}{\left( \frac{n}{3}!\right)^3}
\end{equation*}
Again, by Stirling's approximation 
\begin{equation*}
\frac{1}{2^{2n}}\binom{2n}{n}\sim \frac{1}{\sqrt{\pi n}}
\end{equation*}
and 
\begin{equation*}
\frac{1}{3^n}\frac{n!}{\left( \frac{n}{3}!\right)^3}=\frac{1}{3^n}\frac{\sqrt{2\pi n}\left(\frac{n}{e}\right)^n}{\left( \sqrt{2\pi \frac{n}{3}}\left(\frac{n}{3e}\right)^\frac{n}{3} \right)^3}=\frac{3\sqrt 3}{2\pi n}
\end{equation*}
so that
\begin{equation*}
u^{(3)}_{2n}< \frac{3\sqrt 3}{2\pi \sqrt \pi}\cdot \frac{1}{n^{3/2}}.
\end{equation*}
The last inequality implies that
\begin{equation*}
u^{(3)}:=\sum_{n=0}^{+\infty}u^{(3)}_{2n}< \frac{3\sqrt 3}{2\pi \sqrt \pi}\sum_{n=0}^{+\infty} \frac{1}{n^{3/2}}<+\infty
\end{equation*}
This is enough to state that the particle will not pass infinitely often through the origin, and thus through every possible point. There is only a less than one probability that the particle will soon or later come back to the origin. The three-dimensional symmetric random walk is transient, so the particle may never return to the origin:
\begin{equation*}
{\cal P}^{(3)}_{0}={\cal P}\left( S_{\tau}\ {\rm ever\ returns\ to\ } (0,0,0)|S_{0}=(0,0,0)\right)=\frac{u^{(3)}-1}{u^{(3)}}<1
\end{equation*}
Alternatively, we can consider two particles performing independent symmetric random walks, the steps occurring simultaneously. Will they ever meet? Let us consider the metric in which the distance between two possible positions is the smallest number of steps leading from one position to the other. If the two particles move one step each, their mutual distance either remains the same or changes by two units, and so their distance either is even at all times or else is always odd. In the second case the particles can never occupy the same position. In the first case it is readily seen that the probability of their meeting at the $n$-th step equals the probability of the first particle's reaching in $2n$ steps the initial position of the second particle. Hence our proposition states that in one and two dimensions but not in three dimensions the two particles are sure infinitely often to occupy the same position. If the initial distance of the two particles is odd, a similar argument shows that they will infinitely often occupy neighboring positions. In one and two dimensions the two particles are certain to meet infinitely often but in three dimensions there is a positive probability that they never meet. It has been showed that ${\cal P}^{(3)}_{0}=0,3405373$ and that this probability decreases with increasing dimensions over the third one \cite{McCrea}.


\begin{thebibliography}{1}
\bibitem{Durrett} R. Durrett, Probability: Theory and Examples, Edition $4.1$, $4$th edition, $2013$, Cambridge University Press

\bibitem{Feller} W. Feller, An Introduction to Probability Theory and its Applications, John Wiley \& Sons, II Edition, $1967$.

\bibitem{Kochetkov} Y. Kochetkov, An easy proof of the Polya's Theorem on Random Walks. ArXiv $1803.00811v1$.

\bibitem{Lange} K. Lange, Polya’s random walks theorem revisited. Amer. Math. Monthly, $122(10)$, $2015$, $1005-1007$.

\bibitem{McCrea} McCrea, W. H. and Whipple, F. J. W.   Random Paths in Two and Three Dimensions. Proc. Roy. Soc. Edinburgh, $60$, $281-298$, $1940$.

\bibitem{Novak} J. Novak, Polya’s random walks theorem. Amer. Math. Monthly, $121(8)$, $2014$.

\bibitem{Steele} J. M. Steele, Stochastic Calculus and Financial Applications, Springer, $2000$.

\bibitem{Williams} D. Williams, Probability with Martingales, Cambridge University Press, $1991$.

\end{thebibliography}
\end{document}